\documentclass[a4paper,8pt,fleqn]{article}
\usepackage[a4paper,top=0.7in,bottom=0.8in,left=.8in,right=.8in]{geometry}
\usepackage{color}
\usepackage{graphicx}

\newcommand{\bl}[1]{#1}                      
\newcommand{\rev}{}
\usepackage{caption}
\definecolor{egyptianblue}{rgb}{0.06, 0.2, 0.65}
\usepackage[colorlinks=true,allcolors=egyptianblue,bookmarks=false,hypertexnames=true]{hyperref} 
\usepackage{newtxtext}
\usepackage{booktabs}
\usepackage{multirow}
\usepackage{fontawesome} 

\usepackage{amsthm}
\usepackage{adforn}
\usepackage[libertine]{newtxmath}
\usepackage[cal=boondoxo,bb=boondox,frak=boondox]{mathalfa}
\usepackage[scaled=.95,type1]{cabin}
\newtheorem{proposition}{Proposition}[section]
\newtheorem{assumption}{Assumption}[section]

\newtheorem{definition}{Definition}[section]
\newtheorem{theorem}{Theorem}[section]
\newtheorem{corollary}{Corollary}[section]

\newcommand\indep{\protect\mathpalette{\protect\independenT}{\perp}}
\def\independenT#1#2{\mathrel{\rlap{$#1#2$}\mkern2mu{#1#2}}}

\usepackage{natbib}
\usepackage{fancyhdr}
\usepackage{titlesec}
\titleformat{\section}[hang]{\large\normalfont\scshape\filcenter}{\thesection.}{1em}{}
\titleformat{\subsection}[hang]{\normalfont\itshape\filcenter}{\thesubsection.}{1em}{}
\titleformat{\paragraph}[runin]{\itshape}{\theparagraph.}{1em}{}

\usepackage{authblk}
\begin{document}
 \title{\bfseries Functional independent component analysis by choice of norm: \\ a framework for near-perfect classification}
   \author[1,2,3]{Marc Vidal\thanks{ \faEnvelope \ \href{mailto:marc.vidalbadia@ugent.be}{marc.vidalbadia@ugent.be} (M.V.)  Dept of Data Analysis, Henri Dunantlaan 1, Ghent, Belgium \\ \textbf{Funding}: This research was partially supported by the Methusalem funding from the Flemish Government and the project FQM-307 of the Government of Andalusia (Spain). We also acknowledge the financial support of Agencia Estatal de Investigación, Ministerio de Ciencia e Innovación (grant number: PID2020-113961GB-I00)
and the IMAG María de Maeztu grant CEX2020-001105-M/AEI/10.13039/501100011033. \\

\noindent \textbf{This is a preprint version}. This is a revised version of the preprint originally posted on arXiv and later published in:
\href{https://doi.org/10.1007/s11634-024-00622-5}{https://doi.org/10.1007/s11634-024-00622-5}.} }
   \author[1]{Marc Leman}
   \author[2]{Ana M. Aguilera}
   \affil[1]{\scshape Ghent University, Ghent, Belgium }
   \affil[2]{\scshape University of Granada, Granada, Spain}
   \affil[3]{\scshape Max Planck Institute, Leipzig, Germany}
   \date{}

\maketitle

\begin{abstract}
We develop a theory for functional independent component analysis in an infinite-dimensional framework using Sobolev spaces that accommodate smoother functions. The notion of penalized kurtosis is introduced motivated by Silverman's method for smoothing principal components. This approach allows for a classical definition of independent components obtained via projection onto the eigenfunctions of a smoothed kurtosis operator mapping a whitened functional random variable. We discuss the theoretical properties of this operator in relation to a generalized Fisher discriminant function and the relationship it entails with the Feldman-Hájek dichotomy for Gaussian measures, both of which are critical to the principles of functional classification. The proposed estimators are a particularly competitive alternative in binary classification of functional data and can eventually achieve the so-called near-perfect classification, which is a genuine phenomenon of high-dimensional data. Our methods are illustrated through simulations, various real datasets, and used to model electroencephalographic biomarkers for the diagnosis of depressive disorder.
\end{abstract}
\noindent 
{\it Keywords:} Depressive disorder, EEG, Functional classification, Feldman-H\'{a}jek dichotomy, Kurtosis, Picard condition, Whitening operator

\section{Introduction}\label{sec1}
Kurtosis extremization is a common problem addressed in the context of independent component analysis (ICA), a broadly used reduction technique that assumes a random vector is a linear mixture of components ``as independent as possible”. Compared to principal component analysis (PCA),  ICA goes beyond linear decorrelation by removing higher-order dependencies via prior reparametrization to unit covariance. This pre-processing step, known as \textit{whitening},  standardizes the variance of {the observed data} to remove scale ambiguity and {enhance} statistical independence. Since the application of ICA is often considered to be pointless under certain normality assumptions, one aims to find in the empirical data, say, interesting non-Gaussian projections. Enhancing non-Gaussianity (e.g., as measured by kurtosis) is therefore considered a pathway to achieve statistical independence, a principle supported by the central limit theorem and other precepts from information geometry \citep{Cardoso22}.

Here, we work in the context of functional data analysis, where random objects (functions, images, shapes...) are inherently assumed infinite-dimensional; see \cite{Wang16,Staicu23} for comprehensive reviews on the subject. Despite in the past decades ICA has spawned notorious interest in the machine learning and statistical literature, its functional counterpart has, in contrast, received relatively limited attention. Methods based on reproducing kernel Hilbert spaces (RKHS) were among the earliest to exploit ICA within an infinite-dimensional feature space \citep{Bach02}. However, it was not until the early 2010s that a theoretical framework for ICA within the context of infinite-dimensional Hilbert spaces was first presented in a conference proceedings \citep{Gutch12}. Since then, few works have appeared. As most of these studies have established, the ICA estimation process in functional settings can be summarized in the following steps: $(i)$ whitening the functional data, $(ii)$ rotation via scatter operators and spectral decomposition, $(iii)$ projection onto the space generated by the operator's eigenfunctions and expansion. Typically, this procedure is described through a demixing and mixing transformation, the latter mapping the $H$-valued functional ICs into the original space of functions. Functional ICA is frequently applied in areas such as outlier detection, anomaly identification, and the classification or clustering of functional data \citep{Pena14,Li16,Vidal21}.

\subsection{ICA for classification of functional data}
Extending ICA to infinite-dimensional Hilbert spaces presents additional challenges due to the compactness of the covariance operator, resulting in its inverse being unbounded and allowing only partial standardization. Functional ICA has been then defined through truncated Karhunen–Loève (KL) type representations, or penalized versions of it as a natural method of regularization and control of the smoothness \citep{Li16,Vidal21}. Restricting the estimation of ICs to the subspace spanned by the first functional principal factors can be useful for enhancing the estimation of non-Gaussian ICs: according to the conditions for IC separability \citep{Theis04,Gutch12}, limiting the presence of Gaussian sources may lead to optimal estimation of non-Gaussian ICs. Nevertheless, here we argue that ICA should be used in the opposite direction when addressing functional classification problems. We show that the way to achieve high classification accuracy relies both in ``Gaussianizing” the data as in exploiting the infinite-dimensionality of the space, potentially through nonlinear approaches such as those suggested in this paper. For a comprehensive and up-to-date review of functional classification, refer to \cite{Wang24}; see also \cite{Baillo11}.

A critical step in functional ICA is estimating the kurtosis operator \citep{Pena14}, whose theoretical properties have been studied in the context of a mixture of two Gaussian random processes with same covariance operator but differing in mean structure (for related work, see \cite{PENA2000257, Pena10}). The authors provided a proof that the spectral decomposition of this operator provides an estimator of a generalized Fisher’s linear discriminant subspace. Here, we aim to revisit and further elaborate on these results in relation to what has been described as asymptotic perfect classification \citep{Delaigle12}. Achieving near-perfect classification is possible when the sample size increases and the data are projected onto a suitable one-dimensional function space, such that the probability of misclassification of new observations converges to zero. It has been proven that the  mechanisms underlying this phenomenon can be explained by the Feldman-H\'{a}jek  dichotomy for Gaussian measures \citep{Berrendero18}. We establish a link between these results and the kurtosis operator, which is deduced from the dual properties of this operator. \bl{For mixtures of two homoscedastic Gaussian random elements with equal prior probabilities, we show that the kurtosis of the projection scores tends to 1 exactly when asymptotically perfect classification becomes attainable}, specifically along solution paths that extend beyond Picard’s condition \citep{Nashed74}. The implications of the current findings are noteworthy, as they offer a potential way to significantly simplify classification problems.

\subsection{Relevance of the proposed method}
We present a functional ICA that extends prior formulations to more general spaces through the use of Sobolev Hilbert space norms and inner products. Silverman's method \citep{Silverman96} is adapted by introducing a roughness penalty into the orthonormality constraint of the eigenfunctions of a kurtosis operator mapping a whitened random variable that is not necessarily confined in a finite-dimensional space. The proposed penalized approach {can preserve} the orthonormality property on the estimated projection scores regardless of the degree of smoothing used, which follows from a result in \cite{Ocana99, Vidal21}. While it is well known that Tikhonov-type regularization offers certain advantages for addressing inverse problems, our approach recasts on differential regularization to preserve the orthonormality property typical of ICA procedures. Furthermore, \cite{Hosseini14} showed that incorporating a Tikhonov term into Silverman's method can enhance the estimation of higher-order eigenfunctions but notably worsens the estimation of lower-order ones, compared to both non-smoothed and Silverman’s smoothed estimators. \bl{Differential regularization does not truncate the spectrum, but reshapes the effective spectral problem through smoothness-dependent attenuation}, which can indirectly mitigate issues related to near-zero eigenvalues by promoting stability and smoothness.

In the same spirit as in \cite{Delaigle12}, we advocated for the centroid classifier integrating not only functional PCA estimators but a variety of functional ICA estimators based on kurtosis. Since we focus on binary classification, in a setting where the two populations differ only in terms of means, the centroid classifier is optimal in the range of functional data, very simple to apply and interpretable. Ultimately, our approach is not driven by minimizing classification error, but rather by the probabilistic principles derived from the properties of kurtosis. We further extend the simulation results in \cite{Pena14} to show the well behavior of our estimators. In view of considerable levels of noise or hypervariability, all tested classifiers tend to lose their properties, so we alternatively propose to use other possibilities arising from smoothing the data per group. By applying this form of nonlinear smoothing, we demonstrate that it is possible reach high accuracy in classification problems. Finally, with the analysis of a neuroscientific dataset, we show that our methods are an effective avenue for modeling electroencephalographic biomarkers of neuropsychiatric disorders. 


\section{Kurtosis based functional ICA} \label{sec2}

\subsection{Independence in infinite dimensions}
We next assume a common scenario in functional data where $H$ is a separable space of real-valued functions on a closed interval $\mathcal{I}=[0,T], T>0$, with inner product operator $\langle\cdot,\cdot\rangle: H \times H \rightarrow \mathbb{R}$ and norm $\|\cdot\|: H \rightarrow [0,\infty)$, therefore a space with Hilbert structure. Given a basic probability space  $(\Omega, \mathcal{A}, \mathbb{P})$, a $H$-valued random variable is the mapping $X:\Omega \rightarrow H$ that is  $\mathbb{B}_H$-measurable, where $\mathbb{B}_H$ is the $\sigma$-field generated by the class of all open subsets of $H$.

While independence is a well-established concept to conventional probability, presumably a generalized theory is yet rather scarce. Hence, before delving into the model specifics, we first seek to build preliminary intuition on the notion of independence in Hilbert spaces. For a classical definition regarding finite sequences of $H$-valued random variables, we refer to \cite{Laha79}, Definition 7.13. However, current definition is nuanced when considering the infinite-dimensional nature of a $H$-valued random variable. The following definition is adapted from \cite{Gutch12}.

\vspace*{10pt}
\begin{definition} \label{defindp} 
Let $ (e_j)_{j\in\mathbb{N}}$ be an orthonormal basis  in $H$. We say that a random variable $X:\Omega \rightarrow H$ is independently separable with respect to $ (e_j)_{j\in\mathbb{N}}$, if for any finite index $J=\{j_1,\dots,j_M\}$, the following assumptions hold: 
\begin{enumerate}
\item $\left\langle X, e_{j_k}\right\rangle \indep\left\langle X, e_{j_{\ell}}\right\rangle,$ for all $k \neq \ell \ (k, \ell \in J)$\footnote{It is worth mentioning that the symbol $\indep$, denoting stochastic independence, can be credited to Gustav Elfving, a Finnish statistician, probabilist, and mathematician, who first introduced it in some lecture notes around 1949–1950; see \cite{Nordstrom19}.},
\item $\mathsf{P}_{ e_{J}}(X) \indep \mathsf{P}_{ e_{\overline{J}}}(X),$
where $\mathsf{P}_{ e_{\overline{J}}}(X)=(I_H - \mathsf{P}_{ e_J})(X)$ with $I_H$ the identity operator on $H$ and $ \mathsf{P}_{ e_J}(X)=\sum_{j\in J} \langle X, e_j\rangle e_j$. 
\end{enumerate}
\end{definition} 
\noindent  Assumption 2 in Definition \ref{defindp} reads as ``the orthogonal projection to the subspace spanned by a subset of the basis is independent to the projection to the complement of the subspace” and, while cumbersome to prove analytically, it bears critical implications in this study. A desirable property of a projection is that  $\mathsf{P}_{ e_{J}}(X) \perp \mathsf{P}_{ e_{\overline{J}}}(X),$ which holds for orthonormal basis even when the projection is truncated, as both components remain orthogonal by construction. However, Assumption 2 goes beyond linear independence and rather entails that both $\mathsf{P}_{ e_{J}}(X)$ and $ \mathsf{P}_{ e_{\overline{J}}}(X)$ are mutually independent to define $X$ as independent, connecting again to Laha and Rohatgi's definition.  This constraint was relaxed in \cite{Li16,Vidal21} by converting factorizations of the Karhuenen-Loève expansion into independent subspaces of finite dimension, the latter exploiting an infinite number of subspaces  by the use of smoothing. Let us further note that if assumption 1 in Definition \ref{defindp} holds, then $(S_{j_k})_{k=1}^M$, where $S_k=\langle X, e_{k}\rangle e_{k} \ (k \in J),$ are mutually independent $H$-valued functional components. This notion of independence will facilitate our understanding of the functional ICA model proposed in the following sections.

\subsection{RKHS structure and whitening}
Consider the Banach space of all bounded operators mapping $H$ into itself, denoted by  $\mathcal{B}_H$. 
If $\mathsf{A} \in \mathcal{B}_H$ maps bounded subsets of $H$ to relatively compact subsets and has a discrete spectrum that accumulates at zero, then $\mathsf{A}$ is compact. A compact operator is of trace-class if satisfies  the trace condition $\operatorname{tr}(|\mathsf{A}|) = \sum_{ j = 1}^{\infty}  \sigma_j(\mathsf{A}) < \infty $, where $\sigma_j(\mathsf{A})$ are the singular values of the operator. A trace-class operator is Hilbert-Schmidt if  $\operatorname{tr}(\mathsf{A}\mathsf{A}^*)=\sum_{ j = 1}^{\infty}\left\|\mathsf{A}  e_{j}\right\|^{2} <\infty$, for any orthonormal basis $ (e_j)_{j\in\mathbb{N}}$ of $H$. The space of these operators is denoted by $\mathcal{B}_{2,H}$, and has associated norm $\left\|\mathsf{A} \right\|_{\mathrm{HS}}=\operatorname{tr}(\mathsf{A}\mathsf{A}^*)^{1/2}$. Given two functions $f,g \in H$, if we define the operator  $(f \otimes g )$ by $(f \otimes g)(x)=\langle x, f\rangle g$ for all $x\in H$, the tensor product $\otimes$ induces a continuous, linear operator of rank 1 and therefore a Hilbert-Schmidt operator. In the subsequent sections, we will make use of a property of tensor products:  $(h \otimes g )(g \otimes f )=\|g\|^2(h \otimes f),$ where  $h\in H$.

\vspace*{10pt}
\begin{assumption} \label{A0}
$\mathbb{E}\|X\|^{2} <\infty$ ($X$ has finite second order moments).
\end{assumption}
\vspace*{10pt}
Under Assumption \ref{A0} consider a $\mathbb{B}$-measurable mapping $\Gamma\in \mathcal{B}_{2,H}$ defined by its action $X \mapsto \Gamma_X$,  where $\Gamma_X=\mathbb{E}(X\otimes X)$ and $X$ is assumed with zero mean. The operator $\Gamma_X$, called the covariance operator of $X$, is a positive, self-adjoint and trace-class operator admitting the spectral representation
\begin{equation*}
\Gamma_X=\sum^\infty_{j=1} \lambda_j (\gamma_j \otimes \gamma_j)=\sum^\infty_{j=1} \lambda_j \mathsf{P}_{\gamma_j},
\end{equation*}
where $(\lambda_j)_{j\geq 1}$ is a sequence of positive eigenvalues in decreasing order, $(\gamma_j)_{j\geq 1}$ their associated eigenfunctions and $(\mathsf{P}_{\gamma_j})_{j\geq 1}$ a sequence of projectors onto the space generated by each $\gamma_j$ \citep{Laha79}. Further, we make the following assumption:
\vspace*{10pt}
\begin{assumption} \label{A1}
$\operatorname{ker}(\Gamma_X)=\{0\}$ (the kernel of $\Gamma_X$ is null).
\end{assumption}
\vspace*{10pt}
\noindent Assumption \ref{A1} implies that $\Gamma_X$ is strictly positive, so that all $\lambda_j > 0$. \bl{The injectivity of $\Gamma_X$ allows one to define the generally unbounded inverse $\Gamma_X^{-1}$ on $\mathrm{ran}(\Gamma_X)$. We write $\Gamma_X^{\dagger}$ and $\Gamma_X^{1/2\dagger}$ for the corresponding Moore--Penrose inverses on their natural domains; throughout, $\dagger$ denotes the Moore--Penrose inverse of an operator.} Since $\Gamma_X$ is a positive operator, we might also consider the square root of $\Gamma_X$, defined by
\begin{equation*} 
\Gamma_X^{1/2}=\sum^\infty_{j=1} \lambda_j^{1/2} \mathsf{P}_{\gamma_j}.
\end{equation*}
We know that $\Gamma_X^{1/2}\Gamma_X^{1/2}=\Gamma_X$, which confirms $\Gamma_X^{1/2}$ is Hilbert-Schmidt but not necessarily trace-class (as assumed here).  This operator will play a critical role from here on. Since the aim of functional ICA is to enhance independence via higher-moment operators, it is of interest to free from dependencies and scale $X$. The way to do that is by transforming $X$ into a suitable random variable whose covariance operator is the identity. This is possible through a whitening transformation, which can be generally defined via the factorization  $\Gamma_X^{\dagger}=\Psi \Psi^*$, with $\Psi = \Gamma_X^{1/2 \dagger}$ as the natural choice, providing a basis for establishing conditions for its existence. For $f \in H$, the necessary and sufficient condition for \bl{$f$ to belong to $\mathrm{ran}(\Gamma_X^{1/2})$} is that
\begin{equation} \label{cond1}
    \sum_{j=1}^{\infty}\left|\left\langle f,\gamma_{j}\right\rangle \right|^{2}/\lambda_{j}<\infty,
\end{equation}
which establishes that the coefficients $\left|\left\langle f,\gamma_{j}\right\rangle \right|^{2}$ must decay faster than the corresponding singular values $\lambda_j$ to ensure convergence (see Proposition 2.2 in \cite{Nashed74}). This condition is typically formulated with the covariance operator’s eigenvalues squared, making it more stringent. The requirement that the RKHS norm obtained from \eqref{cond1} must be finite is usually known as the Picard condition. Note that this norm is induced by the inner product
\begin{equation} \label{innprod}
\langle f, g\rangle_{H(\Gamma_X)}=\sum_{j=1}^{\infty} \lambda_j^{-1}\langle f, \gamma_j\rangle \langle g, \gamma_j\rangle=\langle\Gamma_X^{1 / 2 \dagger}f, \Gamma_X^{1 / 2 \dagger}g\rangle,
\end{equation}
where \bl{$H(\Gamma_X)=\mathrm{ran}(\Gamma_X^{1/2})$ denotes the RKHS of $\Gamma_X$}. The Picard condition is critical in functional settings, as it underpins the definition of canonical correlation, the Mahalanobis distance, whitening transformations  (therefore ICA) in an infinite-dimensional framework \citep{Berrendero20,Vidal22}.
 
In the following sections, we use the notation $\mathbb{X}=\Psi\left\{\mathsf{P}_\gamma(X)\right\}$, where $\mathsf{P}_\gamma$ is the projection operator onto the span of the $\gamma$'s. We denote the space of any whitening counterpart $\mathbb{X}$ of $X$ by $\mathbb{M}$, which is assumed to be equipped with the standard inner product. The identification of optimal whitening transformations for functional data has been recently studied in \cite{Vidal22}. The whitening operators defined there are tested in our numerical simulations and real data examples.

Next, we present the functional ICA model, and we introduce the concept of penalized kurtosis, which addresses the main issue of defining functional ICA within an infinite-dimensional framework.

\subsection{Model outline}
 The aim of functional ICA is to enhance the estimation of independent components via orthogonal rotations of $\mathbb{X}$. Let $U_\mathbb{M}$ denote the class of all unitary operators in $\mathcal{B}_\mathbb{M}$.  The functional IC model can be expressed as
\begin{equation}
\label{MFICA}
    \mathcal{W}(X)=\mathcal{U}\Psi(\mathsf{P}_{\gamma} X)=Z,
\end{equation}
where $\mathcal{W}$ is commonly known as the demixing operator, $\mathcal{U}\in U_\mathbb{M}$, and $Z$ is a $H$-valued element with independent component functions  satisfiying $\Gamma_Z=\mathsf{P}_{\overline{\operatorname{ran}}(\mathcal{U}\Psi)}$. Another way to  see model \eqref{MFICA} is by means of a mixing operator \bl{$\mathsf{A}$} which corresponds to
\begin{equation} \label{icamodel2}
    X=\mathsf{A}(Z),
\end{equation}
where \bl{$\mathsf{A}=\Psi^{-1}\mathcal{U}^{*}$}, and naturally, $\mathsf{A}^\dagger=\mathcal{W}$. Briefly, we now discuss some facts about the nature of the mixing operator $\mathsf{A}$. It turns out that this operator is severely unidentified and ill-conditioned. \cite{Gutch12} showed that if none of the components of $Z$ are Gaussian but independent, then $\mathsf{A}$ maps each component of $Z$ to a single component of $X$, \bl{so that $\mathsf{A}$ reduces to a permutation of the components, up to sign and scale}. Otherwise, if there exist two indices $i \neq j \in \mathbb{N}$ such that $\langle e_i, \mathsf{A} e_k\rangle \neq 0 \neq\langle e_j, \mathsf{A} e_k\rangle$ (this operation describes specific non-orthogonal relationships with regard to an arbitrary orthonormal basis of $H$), then, the components of $Z$ must be Gaussian, and any mixing of the components of $Z$ into more than one component of $X$ will be Gaussian too (because we assume independence). This suggests that reaching independence is rather restrictive, and in the end, the aim of functional ICA is to find strategies to estimate components ``as independent as possible”.

A key question in functional ICA is therefore how to determine $\mathcal{U}$. One approach to obtain this operator is through the eigendecomposition of a kurtosis operator. Under the assumption $\mathbb{E}\|X\|^4<\infty$, we can establish the existence of a mapping $\mathsf{K} \in \mathcal{B}_{2,H}$ with the action $X \mapsto \mathsf{K}_X$, defined as 
\begin{equation} \label{eq:stkurt}
\mathsf{K}_X=\mathbb{E}\left[(X\otimes X)^2\right].
\end{equation}
We refer to \eqref{eq:stkurt} as the kurtosis operator of $X$ and, as showed in \cite{Li16}, $\mathsf{K}_X$ is self-adjoint, positive definite, trace-class and unitary equivariant with respect to $U_H$ (the class of all unitary operators in $\mathcal{B}_H$). By the tensor product properties, $\mathsf{K}_X$ is often conveniently expressed as the weighted operator $\mathsf{K}_X=\mathbb{E}\left[\|X\|^2(X\otimes X)\right]$.
 
By analogy to the multivariate case, the independent components could be defined as $\langle \mathbb{X},\psi_j\rangle$, where
 $(\psi_{j})_{j\geq 1}$ is an orthonormal family called independent component weight functions obtained  by solving $\psi_j=$ argmax$_f\mathtt{kurt}(\langle \mathbb{X},f\rangle)$ and subject to $\|f\|^{2}=1,\langle f,\psi_{j}\rangle=0$. This way, the kurtosis based functional ICA is determined by the solutions to the eigenproblem
\begin{equation*}
\mathsf{K}_{\mathbb{X}}\left(\psi_{j}\right)=\kappa_{j} \psi_{j}.
\end{equation*}
The independent component scores $\xi_j=\langle \mathbb{X}, \psi_{j}\rangle$ are then generalized linear combinations of $\mathbb{X}$ with maximum kurtosis satisfying $\mathtt{kurt}(\xi_j)=\langle \mathsf{K}_{\mathbb{X}}\psi_j,\psi_j\rangle = \kappa_j$, where $\mathtt{kurt}$ is understood in terms of how $\mathsf{K}_{\mathbb{X}}$ is defined (note that this does not coincide with the classical univariate coefficient of kurtosis).
However, 
note that
\begin{equation*}
\mathbb{E}\|\mathbb{X}\|^2=\mathbb{E}\left(\sum_{j=1}^{\infty} \frac{\left|\left\langle X, \gamma_j\right\rangle\right|^2}{\lambda_j}\right)=\sum_{j=1}^{\infty} \frac{1}{\lambda_j}\left\langle\gamma_j, \mathbb{E}\left(\left\langle X, \gamma_j\right\rangle X\right)\right\rangle=\sum_{j=1}^{\infty} \frac{\lambda_j}{\lambda_j}\left\langle\gamma_j, \gamma_j\right\rangle=\infty,
\end{equation*}
($\mathbb{X}$ has infinite variance). \bl{In fact, assuming $\mathbb{X}$ exists,  $\Gamma_\mathbb{X}=\mathsf{P}_{\overline{\operatorname{ran}}(\Psi)}$}, where the projection operator \bl{$\mathsf{P}_{\overline{\operatorname{ran}}(\Psi)}$} defines the identity in the space; see Sect. 2 in \cite{Vidal22}. \bl{Thus, on $\overline{\operatorname{ran}}(\Psi)$, $\Gamma_{\mathbb{X}}$ acts as the identity operator; consequently, if this space is infinite-dimensional, $\operatorname{tr}(\Gamma_{\mathbb{X}})=\infty$.} 

The primary challenge in defining ICA within an infinite-dimensional framework is that higher order moments require the integration of lower-order terms. \bl{Therefore, full whitening in infinite dimension does not yield a square-integrable $H$-valued element with identity covariance, and the unregularized kurtosis operator of $\mathbb{X}$ cannot in general be treated as a trace-class fourth-moment operator, which would require $\mathbb{E}\|\mathbb{X}\|^4<\infty$.} Moreover, it is unlikely that the paths of $X$ lie within the RKHS associated with their covariance kernel: the probability of the sample paths belonging to the RKHS is often zero due to insufficient regularity \citep[Theorem 7.3]{Milan01}. Overall, these issues suggest that functional ICA must be defined through a truncated KL expansion or by the use of some kind of regularization in the norm generated by \eqref{cond1}  (in case one would like to minimize the effects of such truncation, as argued in  \cite{Berrendero20}). In the next subsection, we will demonstrate that it is possible to define the kurtosis operator of a smoothed version of $\mathbb{X}$, allowing a definition of the ICA model similar to that in \eqref{icamodel2}, i.e., in an infinite-dimensional framework.

\subsection{Penalized kurtosis}
Here, we extend the method introduced in \cite{Silverman96} for smoothing principal component estimates to functional ICA by introducing the notion of penalized kurtosis. In doing so, we aim at controlling the smoothness of the IC weight functions by incorporating roughness penalties through a linear differential operator. \bl{In the empirical setting,} the small perturbations produced by this type of penalty will presumably have a regularizing effect on the lower order eigenelements by adjusting possible distortions.

In what follows, we shall consider that $\mathbb{M}_\theta$ is a closed subspace of continuously differentiable functions with weighted Sobolev inner product and corresponding norm
\begin{equation*}
\left\langle f,g\right\rangle _{\theta}=\left\langle f,g\right\rangle+\theta\left\langle \mathsf{D}_r f,\mathsf{D}_r g\right\rangle,\quad \|f\|^2_\theta=\left\langle f,f\right\rangle _{\theta},
\end{equation*}
where $\theta \in \mathbb{R}_+$ is a penalty parameter and \bl{$\mathsf{D}_r$ a densely defined closed differential operator of order $r$ with the appropriate domain}.

In order to obtain the main results presented in the next section, let us consider Proposition 4.3 in \cite{Ocana99} for the covariance operator (functional PCA). Assuming continuity of the usual inner product $\left\langle \cdot,\cdot\right\rangle$  in terms of the modified inner product $\left\langle \cdot,\cdot\right\rangle _{\theta}$,  there exists an  injective, symmetric and  bounded operator $\mathsf{S}^2$ such that $\left\langle f,g\right\rangle=\langle \mathsf{S}^{2}(f),g\rangle _{\theta}$. Given these properties, we can assume the existence of a square root operator $\mathsf{S}$, which is also a symmetric, positive definite, satisfying $\mathsf{S}^2=\mathsf{S} \mathsf{S}$, so that
 \begin{equation} \label{rel}
			\left\langle f,g\right\rangle=\langle \mathsf{S}^{2}(f),g\rangle _{\theta} = \left\langle \mathsf{S}(f),\mathsf{S}(g)\right\rangle _{\theta}.
   \end{equation}

{\rev
Let $\mathsf{K}_\theta$ be the operator determined on $\mathbb{M}_\theta$ by
\begin{equation*}
\langle f,\mathsf{K}_\theta g\rangle=\mathbb{E}\left\{\|\mathsf{S}(\mathbb{X})\|^{2}\langle f,\mathbb{X}\rangle\langle g,\mathbb{X}\rangle\right\}.
\end{equation*}
}
Analogously to Silverman's method, the novel penalized ICA approach maximizes
\begin{equation} \label{kurtpen}
\frac{\bl{\langle f,\mathsf{K}_{\theta}f\rangle}}{\left\langle f,f\right\rangle +\theta\left\langle \mathsf{D}_rf,\mathsf{D}_rf\right\rangle }=\frac{\bl{\langle f,\mathsf{K}_{\theta}f\rangle}}{\left\Vert f\right\Vert _{\theta}^{2}},
\end{equation}
for all $f \in \mathbb{M}_\theta, f\neq0$. Note that $\theta$ controls the  roughness of the function $f$ as measured by the penalty $\left\langle \mathsf{D}_r f,\mathsf{D}_rf\right\rangle$. Consequently, one can find a collection of smoothed functions $\psi_{\theta,j}\in \mathbb{M}_\theta$ that maximize \eqref{kurtpen} which is equivalent to solve the following optimization problem:
\begin{equation}
\begin{split} 
\psi_{\theta,1}&={\mathrm{argmax}_{f}}\left\langle f,\bl{\mathsf{K}_{\theta}}f\right\rangle \ \mathrm{s.\ t. \ } \ \|f\|_{\theta}^{2}=1,\\ \label{1PSV}
\psi_{\theta,k}&=
{\mathrm{argmax}_{f}}\left\langle f,\bl{\mathsf{K}_{\theta}}f\right\rangle \ \mathrm{s.\ t. \ } \|f\|_{\theta}^{2}=1,\left\langle f,\psi_{\theta,j}\right\rangle _{\theta}=0,\ \mathrm{for\ all \ } j<k \ \ (k= 2,3,\ldots). 
\end{split}
\end{equation}
The orthonormal condition over the smoothed IC  weight functions is  fixed by the inner product $\langle \cdot,\cdot\rangle _{\theta}$ whereas the kurtosis of the independent components is given by $\langle \cdot,\cdot\rangle $. In this sense, the smoothed IC weight functions form an orthonormal system of $\mathbb{M}_\theta$.

 \vspace*{10pt} 
 \begin{assumption} \label{finite4}
$\mathbb{E}\|\mathsf{S}\mathbb{X}\|^4<\infty.$
 \end{assumption}
 \vspace*{10pt}
 \begin{proposition} \label{Prop01}
{\rev
The smoothed kurtosis operator $\mathsf{K}_{\mathsf{S}(\mathbb{X})}$ has finite trace if and only if $\mathbb{E}\|\mathsf{S}(\mathbb{X})\|^4<\infty$, that is, under Assumption \ref{finite4}.
}
\end{proposition}

\vspace*{10pt}
\begin{proposition} \label{Prop02}
Under Assumption \ref{finite4}, the eigensystem of the smoothed functional ICA denoted by $(\psi_{\theta,j}, \kappa_{\theta,j}) \in \mathbb{M}_\theta \times\mathbb{R},$ is obtained as the solutions to the equation 
\begin{equation*}
\langle f,\bl{\mathsf{K}_{\theta}}\psi_{\theta,j}\rangle =\kappa_{\theta,j}\left\langle \psi_{\theta,j},f\right\rangle _{\theta},
\end{equation*}
which is equivalent to the eigensystem of $\mathsf{K}_{\mathsf{S}^2(\mathbb{X})}$ with the inner product $\langle \cdot,\cdot\rangle_\theta.$ 
\end{proposition}

 \vspace*{10pt}
Equivalently to the results provided in \cite{Ocana99} for the smoothed functional PCA, this algorithm can be regarded as an equivalence between the smoothed functional ICA  and the spectral decomposition of  the kurtosis operator of the half-smoothed process $\mathsf{S}(\mathbb{X})$ with the usual inner product. 

 \vspace*{10pt}
\begin{proposition} \label{Prop0}
$(\psi_\theta,\kappa_\theta) \in \mathbb{M}_\theta \times\mathbb{R}$ is an eigenelement of $\mathsf{K}_{\mathsf{S}^2(\mathbb{X})}$ with $\langle \cdot,\cdot\rangle_\theta$ if and only if $(\mathsf{S}^{-1}(\psi_\theta),\kappa_\theta)$ is an eigenelement of  $\ \mathsf{K}_{\mathsf{S}(\mathbb{X})}$ with $\langle \cdot,\cdot\rangle$.
\end{proposition}

\vspace*{10pt}
\noindent Then, we can establish the following equivalences which follows immediately from Propositions \ref{Prop02} and  \ref{Prop0}.

\vspace*{10pt}
\begin{corollary} \label{Prop1}
The smoothed independent components $(\xi_{\theta,j})_{\ j\in \mathbb{N}_+}$ are equivalently obtained by the following projections:
\begin{enumerate}
    \item $\left\langle \mathbb{X},\psi_{\theta,j} \right\rangle$
    \item $\left\langle \mathsf{S}^2(\mathbb{X}),\psi_{\theta,j} \right\rangle _{\theta}$
    \item $\left\langle \mathsf{S}(\mathbb{X}), \mathsf{S}^{-1} (\psi_{\theta,j}) \right\rangle$.
\end{enumerate}
As a consequence, the orthogonal representation for the half-smoothed whitened variable $\mathsf{S}(\mathbb{X})$ in terms of the independent components is obtained by the expansion
\begin{equation*}
\mathsf{S}(\mathbb{X}) = \sum_{j=1}^\infty \xi_{\theta,j} \mathsf{S}^{-1} (\psi_{\theta,j}).
\end{equation*}
\end{corollary}

\noindent The IC scores $\xi_{\theta,j}$ are then generalized linear combinations of $\mathsf{S}(\mathbb{X})$ with maximum kurtosis satisfying 
 \begin{equation*}
\mathtt{kurt}(\xi_{\theta,j})=\left\langle\mathsf{K}_{\mathsf{S}(\mathbb{X})} \mathsf{S}^{-1}(\psi_{\theta,j}), \mathsf{S}^{-1}(\psi_{\theta,j})\right\rangle=\kappa_{\theta,j}.
 \end{equation*}
Model \eqref{MFICA} is completed by defining \bl{$\mathcal{U}=\sum_{j=1}^{\infty}\left\{\mathsf{S}^{-1}(\psi_{\theta,j}) \otimes e_j\right\}$, where $(e_j)_{j\geq1}$ is a fixed orthonormal basis of $\mathbb{M}$ and the system $\{\mathsf{S}^{-1}(\psi_{\theta,j})\}$ is orthonormal for $\langle\cdot,\cdot\rangle$,} and the operator $ \mathsf{A}$ as
 \begin{equation*}
  \mathsf{A}=(\mathcal{U}\mathsf{S}\Psi +  \mathcal{W}_ \epsilon)^\dagger,
 \end{equation*}
 where $ \mathcal{W}_\epsilon=\mathsf{S}\Psi - \Psi$ is a perturbation of $\mathsf{S}\Psi$.

\section{Theoretical properties}
\subsection{Discriminative properties of the kurtosis operator} \label{sec3.1}
In this section, we examine the theoretical properties of the kurtosis operator as they relate to Fisher's linear discriminant function.  Suppose that  $X$ can be observed as a mixture of two subpopulations $\Pi_k \ (k=0,1)$, which are identified by the binary variable $Y=k$ when $X\in \Pi_k$, and the aim is to assign their sample paths into one of them. Let $\pi_k=\mathbb{P}(Y=k)$ with $\pi_1=1-\pi_0$, and consider that class $k$ has mean function $\mu_k=\mathbb{E}(X | Y=k)$ and  $\mathbb{E}(X) = \mu =\pi_0 \mu_0+(1-\pi_0)\mu_1$. In principle, we do not impose distributional assumptions on $X$, but we consider $\mu_0\neq\mu_1$ and  equal class covariance operators. We will assume the functions in $H$ are square integrable on $\mathcal{I}$. \bl{Throughout Sections \ref{sec3.1}--\ref{SCons}, $(\lambda_j,\gamma_j)_{j\geq1}$ denote the eigenelements of the common within-class covariance operator $\mathcal{R}$.} 

Fisher's discriminant problem consists in estimating a non-zero function $\varphi$ on $\mathcal{I}$ that maximizes the ratio
\begin{equation} \label{DF}
\mathsf{F}(\varphi)=\left\langle \varphi, \Gamma_{X}^{W}\varphi\right\rangle^{-1}\left\langle \varphi, \Gamma_{X}^{B}\varphi\right\rangle,  
\end{equation}
where
$$\Gamma_{X}^{B}=\pi_0(1-\pi_0) \{(\mu_{1}-\mu_{0}) \otimes (\mu_{1}-\mu_{0})\} \quad \mathrm{and} \quad
\Gamma_{X}^{W}=\pi_0\Gamma_{X|Y=0} + (1-\pi_0)\Gamma_{X|Y=1}=\mathcal{R}$$ 
are respectively, the between and within-class covariance operator. The idea in \eqref{DF} is to give large separation to the group means while, at the same time, keeping the variance between groups small. 

We note that, $\mathcal{R}\equiv \sum^\infty_{j=1} \lambda_j \mathsf{P}_{\gamma_j}$ represents the common covariance operator in each population. According to the law of the total covariance, one has that $\Gamma_{X}=\Gamma_{X}^{W}+\Gamma_{X}^{B}$ and then, for a function $f$ with expansion $\bl{f=\sum_{j=1}^{\infty}\mathtt{f}_{j} \gamma_j}$, $\Gamma_{X}$ can be written as
\begin{equation*}
\Gamma_{X}(f)
=\sum_{j=1}^\infty \left(\lambda_j \mathtt{f}_j
+ \pi_0(1-\pi_0)\mathtt{v}_j \sum_{k=1}^\infty \mathtt{v}_{k} \mathtt{f}_{k}\right)\gamma_j,
=\sum_{jk=1}^\infty s_{jk}\mathtt{f}_j\gamma_{k},
\end{equation*}
where $\mathtt{v}_j$ are the coefficients of the mean differences between classes in terms of $\gamma_j$'s (i.e., $\mu_\Delta=\mu_0-\mu_1=\sum_{j=1}^\infty \mathtt{v}_j \gamma_j $) and 
$$
s_{j k}=\left\{\begin{array}{cl}
\lambda_j+\pi_0(1-\pi_0)\mathtt{v}_{j}^2 & j=k \\
 \pi_0(1-\pi_0)\mathtt{v}_j \mathtt{v}_{k} &  j\neq k.
\end{array}\right.
$$

Solutions to (\ref{DF}) are well-known. Here, we briefly debrief them for the sake of clarity.  
\vspace*{10pt}
\begin{proposition}[\cite{Pena14}, Lemma 1] \label{lemmadiscr} 
Suppose that all $\lambda_j>0$. \bl{Whenever the series below defines an element of $H$,} for some constant $c$, 
\begin{equation} \label{eq:discrim}
\varphi=c\sum_{j=1}^\infty\lambda_j^{-1} \mathtt{v}_j\gamma_j
\end{equation}
is the maximizer of $\mathsf{F}$.
\end{proposition}

\vspace*{10pt}
\noindent {\rev In particular, \eqref{eq:discrim} defines an element of $H$ when $\mu_\Delta\in\mathrm{ran}(\mathcal{R})$. For each finite $q$, $\beta^{[q]}=\sum_{j=1}^{q}\lambda_j^{-1}\mathtt{v}_j\gamma_j$ maximizes $\mathsf{F}$ over $\mathrm{span}\{\gamma_1,\ldots,\gamma_q\}$.}

\vspace*{10pt} 
\noindent According to this result, \cite{Pena14} provided a proof that \bl{the kurtosis operator} has a spectral direction equivalent to the function found in Proposition \ref{lemmadiscr}. The proof was given in the context of Gaussian random elements. \bl{The following consequence of the argument in \cite{Pena14} (Theorem 1) will be useful here.}

\vspace*{10pt}
\begin{proposition} \label{Ke} Let $X$ be a mixture of two Gaussian \bl{random elements} with different means and same covariance operator, and let \bl{$\mathbb{X}\equiv \Gamma_X^{1/2\dagger}(X-\mu)$} \bl{be an approximation to a standardized random element}. {\rev Then the argument of \cite{Pena14} implies that the discriminant direction corresponds, after whitening, to an eigenfunction $\psi_\Delta$ of $\mathsf{K}_\mathbb{X}$, with
\begin{equation*}
\psi_\Delta \propto \Gamma_X^{1/2\dagger}\mu_\Delta .
\end{equation*}
Consequently,
\begin{equation*}
\langle X-\mu,\varphi\rangle \propto \langle \mathbb{X},\psi_\Delta\rangle .
\end{equation*}}
\end{proposition}

\subsection{The duality of the kurtosis operator and near-perfect classification} \label{TS4}
Let $(X_{i},Y_{i})$ be independent and identically distributed data pairs drawn from $(X,Y)$, where $Y_i \in \{0,1\} $ is the group label. For $k=0,1$ and  $1 \leq i \leq n_k$, let $X_{ik}$ denote the $i$th function among $X_1, \ldots, X_n$ from the subpopulation $\Pi_k$, where $n_0+n_1=n$ and assume $\pi_0=1/2$ throughout this section. We consider the problem of classifying a newly drawn function, say $X^*$, with the goal of assigning it to one of the two subpopulations, $\Pi_0$ or $\Pi_1$. \cite{Delaigle12} proposes an asymptotic version of the finite sample centroid-classifier with desirable properties we will use in our analyses. The finite version of this classifier assigns $X^*$ to $\Pi_0$ or $\Pi_1$ based on whether the statistic $T_n\left(X^*\right)=\mathtt{dist}^2\left(X^*, \hat{\mu}_1\right)-\mathtt{dist}^2\left(X^*, \hat{\mu}_0\right)$ is positive or negative, respectively. Suppose the above classification rule defined by the distance  $\mathtt{dist}\left(X^*, \hat{\mu}_k\right)=\left|\langle X^*, \beta\rangle-\left\langle\hat{\mu}_k, \beta\right\rangle\right|$, where $\beta$ is a pre-chosen direction defined on $\mathcal{I}$. If we let $n_0,n_1 \to \infty$ and assume $\hat{\mu}_0=0$ and $\hat{\mu}_1\equiv \hat{\mu}_\Delta=\sum_{j=1}^\infty \mathtt{v}_j\gamma_j $, the asymptotic version of $T_n$ is
\begin{equation*}
T^0(X^*)=(\langle X^*, \beta\rangle-\langle \hat{\mu}_\Delta, \beta\rangle)^2-(\langle X^*, \beta\rangle)^2.
\end{equation*}
Delaigle and Hall showed that $T^0$ ensures optimal classification accuracy in the Gaussian setting. Since our data is unlikely to be of such nature, the following theorem addresses the performance of $T^0$ under more general distributional conditions. We discuss its implications concerning Picard's condition and the Feldman-H\'{a}jek dichotomy.
\vspace*{10pt}
\begin{theorem}[\cite{Delaigle12}, Theorem 2] \label{TDel}
Assume $X_{ik} \ (i=1,\dots,n; k=0,1),$ are not Gaussian, $\hat{\mu}_0\neq\hat{\mu}_1$ and $\hat{\Gamma}_{X_{k}}$ are the same.   If $\Pi_0$ and  $\Pi_1$ have prior probabilities $\pi_0$ and $1-\pi_0$ respectively, and $\mu_0=0$, then
\begin{enumerate}
    \item The misclassification probability  for the classifier $T^0$ equals $\mathrm{err}=\pi \mathbb{P}(Q> \nu/2\sigma_Q) + (1-\pi)\mathbb{P}(Q<-\nu/2\sigma_Q)$, where $Q=\langle X-\hat{\mu}_\Delta,\beta\rangle$ and $\nu=\langle \hat{\mu}_\Delta ,\beta\rangle$.
    \item  If $\sum_{j=1}^\infty\lambda^{-1}_j\mathtt{v}^2_j=\infty$, by taking a sequence of classifiers build from $\beta^{[q]}=\sum_{j=1}^{q} \lambda^{-1}_j\mathtt{v}_j\gamma_j$ with $q \to \infty$, the minimal misclassification probability (minimum value of $\mathrm{err})$ tends to $\mathrm{err}_0=0$, and perfect classification is then possible.
\end{enumerate}
\end{theorem}

\noindent Geometrically, from Theorem \ref{TDel} one deduces that asymptotic perfect classification is related to the divergence in norm induced by the inner product defined in \eqref{innprod} as  \bl{$\|\mu_\Delta\|^2_{{H(\mathcal{R})}}$}$=\sum_{j=1}^\infty\lambda^{-1}_j\mathtt{v}^2_j=\infty$, where \bl{$H(\mathcal{R})=\operatorname{ran}(\mathcal{R}^{1/2})$}. Note this is the same as assuming that Picard's condition does not hold (the coefficients $\mathtt{v}^2_j$ decay slower than the corresponding $\lambda_j$s and therefore, \bl{$\mu_\Delta \notin \mathrm{ran}(\mathcal{R}^{1/2})$}. \bl{Observe that, for each finite $q$, the direction $\beta^{[q]}=\sum_{j=1}^{q} \lambda^{-1}_j\mathtt{v}_j\gamma_j$ corresponds, after whitening, to the discriminant eigenfunction described in Proposition \ref{Ke}.}

In the Gaussian setting (see Theorem 1 in \cite{Delaigle12}), the reason of the behavior described in Theorem \ref{TDel} has a probabilistic interpretation by the Feldman-H\'{a}jek dichotomy for Gaussian measures. Two probability measures $m_k, k\in \{0,1\},$ are said to be equivalent ($m_0\sim m_1$) if they are absolutely continuous with respect to one another: i.e., if $m_0(B)=0$ for a Borel set $B \in \mathbb{B}_H$, it holds $m_1(B)=0$  (they have the same zero sets). Conversely, if $m_0(B)=0$ and $m_1(B)=1$, then we say that $m_0$ and $m_1$ are mutually singular $(m_0\perp m_1)$ as $B$ splits in two disjoint sets where $m_0$ and $m_1$ are respectively concentrated. The Feldman-H\'{a}jek dichotomy states that in infinite dimensions, two Gaussian measures have the critical property of being either equivalent or mutually singular. 

\vspace*{10pt}
\begin{theorem}[\cite{daprato14}, Theorem 2.25] \label{HF}
Let $m_k=N(\mu_{m_k},\Gamma_{m_k})$ \\ $(k=0,1)$, be two Gaussian measures on $H$. Then, $m_0\sim m_1$  if and only if, it holds:
\begin{enumerate}
\item  Both measures share the same space, i.e., $\mathrm{ran}(\Gamma^{1/2}_{m})=\mathrm{ran}(\Gamma^{1/2}_{m_0})=\mathrm{ran}(\Gamma^{1/2}_{m_1})$.

\item $\mu_{m_0}-\mu_{m_1} \in \mathrm{ran}(\Gamma^{1/2}_{m})$.

\item $(\Gamma^{-1/2\dagger}_{m_0}\Gamma^{1/2}_{m_1})(\Gamma^{-1/2\dagger}_{m_0}\Gamma^{1/2}_{m_1})^*-I_H$ is a Hilbert-Schmidt operator on $\overline{\mathrm{ran}}(\Gamma^{1/2}_{m})$.
\end{enumerate}
If one of the above conditions is violated, then $m_0\perp m_1$.
\end{theorem}

\vspace*{10pt}
\noindent As a consequence of the Theorem \ref{HF}, we also consider the following result.

\vspace*{10pt}
{\rev
\begin{corollary} \label{cor2}
Suppose the two classes are Gaussian with common covariance $\mathcal{R}$, $\pi_0=1/2$, and $\mu^{[q]}_\Delta\neq0$; write $\mathbb{X}^{[q]}=\Gamma^{1/2\dagger}_{X^{[q]}}(\mathsf{P}^{[q]}_\gamma X-\mu^{[q]})$, and let $\psi_q$ be the eigenfunction of $\mathsf{K}_{\mathbb{X}^{[q]}}$ given by Proposition \ref{Ke}. Put
\begin{equation} \label{pq}
d_q^2=\sum_{j=1}^{q}\frac{\mathtt{v}_j^2}{\lambda_j}, \qquad
p_q=\frac{d_q^2}{4+d_q^2},
\end{equation}
and let $\mathtt{kurt}_q=\kappa_q-(q-1)$ denote the normalized eigenvalue of $\mathsf{K}_{\mathbb{X}^{[q]}}$ at $\psi_q$. Then $\mathtt{kurt}_q=3-2p_q^2$, the remaining normalized eigenvalues are all equal to $3$, the misclassification probability of $T^0$ along $\beta^{[q]}$ is $\mathrm{err}_q=\Phi(-d_q/2)$, where $\Phi$ denotes the standard normal distribution function, and
\begin{equation} \label{gchain}
\mathtt{kurt}_q\to1
\iff d_q^2\to\infty
\iff \mu_\Delta\notin H(\mathcal{R})
\iff m_0\perp m_1
\iff \mathrm{err}_q\to0 .
\end{equation}
In particular $(3-\mathtt{kurt}_q)/2=p_q^2\in[0,1)$.
\end{corollary}

}

\vspace*{10pt}
\noindent Now, recall Theorem \ref{TDel}. Suppose that all $X_i$ are Gaussian via the measures $m_k$. In \cite{Berrendero18} (Theorem 5), it has been proven that \bl{$m_0 \sim m_1 \Leftrightarrow \|\mu_\Delta\|^2_{H(\mathcal{R})}<\infty$ and $m_0\perp m_1 \Leftrightarrow \|\mu_\Delta\|^2_{H(\mathcal{R})}=\infty$}, thus explaining the mechanisms underlying the dichotomy found in \cite{Delaigle12}. This result follows from Theorem \ref{HF} and Parseval's formula as \bl{$\mu_\Delta \in \mathrm{ran}(\mathcal{R}^{1/2})$ if and only if $\|\mu_\Delta\|^2_{H(\mathcal{R})}<\infty$}.

Let us define the sample kurtosis operator
\begin{equation} \label{skurtop}
\hat{\mathsf{K}}_{\mathbb{X}^{[q]}}=\frac{1}{n} \sum_{i=1}^{n}\left\|\mathbb{X}_i^{[q]}\right\|^2\left(\mathbb{X}_i^{[q]} \otimes \mathbb{X}_i^{[q]}\right)
\end{equation}
where $\mathbb{X}_i^{[q]} = \Psi(\mathsf{P}^{[q]}_\gamma X_i)$ with $\mathsf{P}^{[q]}_\gamma$ the rank \bl{$q<n-1$} projection operator onto a set of the span of the $\gamma_j$'s. Next, consider the normalized kurtosis $\mathtt{kurt}(\xi^{[q]}_j)= \hat{\kappa}_j-(q-1),$ where $(\hat{\kappa}_j,\hat{\psi}_j)$ denote the eigenvalues and eigenfunctions of $\hat{\mathsf{K}}_{\mathbb{X}^{[q]}}$.
The term $(q-1)$ serves as a normalization factor, ensuring that the kurtosis equals 3 in the Gaussian case\bl{; this follows from the computation $\kappa_u=q+2$ for directions orthogonal to the discriminant one in the proof of Corollary \ref{cor2}}. \bl{By Corollary \ref{cor2}, the population normalized spectrum is constant at $3$ except in the discriminant direction, where it equals $3-2p_q^2$ and tends to $1$ as $d_q^2\to\infty$}, indicating the presence of bimodality; \bl{the corresponding sample eigenspaces are consistently estimated under the fixed-$q$ conditions of Sect.~\ref{SCons}}. 
{\rev
At the population level, Corollary \ref{cor2} gives
\begin{equation*}
\frac{3-\mathtt{kurt}_q}{2}=p_q^2\in[0,1),
\end{equation*}
and this quantity converges to $1$ precisely when $d_q^2\to\infty$. Consequently, along the population discriminant sequence, $\mathtt{kurt}_q\to1$ if and only if the misclassification probability approaches $\mathrm{err}_0=0$.
} 

Not surprisingly, low kurtosis has been previously associated, with some reservations, to bimodality in symmetric distributions. \bl{In the Gaussian setting, Corollary \ref{cor2} shows that the population kurtosis spectrum determines on which side of the Feldman--H\'{a}jek dichotomy the two class measures lie, and hence the attainable classification error.} The \bl{minimization of the kurtosis} naturally occurs via the eigendecomposition of the kurtosis operator, but other operators can be considered. For instance, solutions to choose the best $\beta$ can lie arbitrarily in the tails of the principal component expansion \citep{Delaigle12}.
In cases of high regularity in $X$, one can expect that solutions concentrate towards the first principal component. A toy example in Appendix B illustrates these facts. We therefore suggest that the kurtosis coefficient can be used to assess the optimality of the $\beta$ estimators, significantly simplifying the computational cost compared to other cross-validation techniques. \bl{It is natural to conjecture that} in the case of more than two populations, the kurtosis coefficient of the projection scores for all pairs of distinct groups should converge to $1$ to reach high classification accuracy. Classification techniques for the current case will be assessed in future studies.

\subsection{Consistency} \label{SCons}
In this section, we show that the eigendecomposition of the sample smoothed kurtosis provides a consistent estimator of the Fisher subspace. There are two challenges in pursuing this: ($i$) the sample covariance operator $\hat{\Gamma}_X$ is non-injective, not even if $\Gamma_X$ is, therefore $\hat{\Gamma}_X$ can have at most $n-1$ non-zero eigenvalues ($ii$) the kurtosis operator $\mathsf{K}_{\mathbb{X}}$ is not trace-class and therefore does not admit a convergent spectral representation in Mercer's sense. In light of these issues, a possible approach to study the consistency of our estimators is by leveraging finite-dimensional projections. Naturally, all good results depend on selecting an appropriate truncation point in the KL expansion. Since Picard's criterion provides a way to assess the ill-posedeness of the precision, if $q$ represents this point, one could chose it as the $q$ for which the decay rates of $\lambda_j$ and $\left|\left\langle X, \gamma_j\right\rangle\right|^2$ become comparable to get closer to the limit of divergence. According to our principles, this will make our classifier optimal. To this end, $q$ can be determined by
\begin{equation*}
q=\arg \min \left\{\sum_{j>q} \log \left(\frac{\left|\left\langle X, \gamma_j\right\rangle\right|^2}{\lambda_j}\right) \text { s.t. } \lim _{j \geq q+\delta} \log \left(\frac{\left|\left\langle X, \gamma_j\right\rangle\right|^2}{\lambda_j}\right)=0\right\},
\end{equation*}
where $\delta \in \mathbb{N}$ allows to improve the numerical stability of the truncation beyond the point the logarithm of the ratio approaches zero, ensuring better control over tail behavior. Then, we can consider the following assumption to validate the consistency of the estimators:
\vspace*{10pt}
\begin{assumption} \label{Awhite}
$\| \hat{\Gamma}_{\mathbb{X}^{[q]}} - I_{H^{[q]}}\|_{\mathrm{HS}}=\mathcal{O}_{\mathbb{P}}\left(\frac{\log(q+\delta)}{n}\right)$.
\end{assumption}
\vspace*{10pt}
\noindent Here, $\hat{\Gamma}_{\mathbb{X}^{[q]}}$ represents the sample covariance of the standardized $q$-dimensional KL expansion of the data, and $\Gamma_{\mathbb{X}^{[q]}} = I_{H^{[q]}}$ is its theoretical counterpart. This assumption is quite standard in functional settings, and it implies that as the sample size $n$ grows, the rate of convergence improves more rapidly. We show that, under mild conditions, these rates are achieved by most   commonly used whitening transformations (see Sect.~\ref{ss5.2}).

Next, consider the kurtosis operator $\mathsf{K}_{\mathbb{X}^{[q]}}=\mathbb{E}\left\{\|\mathbb{X}^{[q]}\|^2\left(\mathbb{X}^{[q]} \otimes \mathbb{X}^{[q]}\right)\right\}$, which is naturally strictly positive definite and trace-class by finite-dimensional dependency. Denote the empirical estimator of the smoothed kurtosis operator a $\hat{\mathsf{K}}_{\mathsf{S}^2(\mathbb{X}^{[q]})}$ (this is the smoothed counterpart of $\hat{\mathsf{K}}_{\mathbb{X}^{[q]}}$ defined in \eqref{skurtop}). The eigensystem of $\hat{\mathsf{K}}_{\mathsf{S}^2(\mathbb{X}^{[q]})}$ is  denoted by $(\hat{\kappa}_{\theta,j},\hat{\psi}_{\theta, j})_{j = 1}^q$. Now, order the eigenvalues of $\mathsf{K}_{\mathbb{X}^{[q]}}$  and $\hat{\mathsf{K}}_{\mathsf{S}^2(\mathbb{X}^{[q]})}$ as a decreasing sequences, so that $\kappa_j \geqslant \kappa_{j+1}$ and $\hat{\kappa}_{\theta,j} \geqslant \hat{\kappa}_{\theta,(j+1)}$. Define \bl{$\mathcal{Y}_j=\left\{1\leqslant l \leqslant q: \kappa_l=\kappa_j\right\}$} and let $\mathcal{y}_j=\left|\mathcal{Y}_j\right|$, where $|\mathcal{Y}|$ represents the cardinality of the finite set $\mathcal{Y}$. Further, for each $\kappa_j$ and $\hat{\kappa}_{\theta,j}$, define the following projector operators:
\begin{equation*}
\bl{\mathsf{P}_{\psi_j}=\sum_{l \in \mathcal{Y}_j}( \psi_l \otimes \psi_l), \quad \mathsf{P}_{\hat{\psi}_{\theta,j}}=\text{orthogonal projection of } H^{[q]} \text{ onto } \operatorname{span}\{\hat{\psi}_{\theta,l}: l \in \mathcal{Y}_j\},} 
\end{equation*}
respectively; \bl{the second is described as a projection rather than as a sum of rank-one terms because the $\hat{\psi}_{\theta,l}$ are orthonormal for $\langle\cdot,\cdot\rangle_\theta$, not for $\langle\cdot,\cdot\rangle$.}
\vspace*{5pt}
\begin{assumption} \label{Asob} The eigenfunctions $\psi_j$ of $\mathsf{K}_{\mathbb{X}^{[q]}}$ satisfy $\int_{\mathcal{I}} \mathsf{D}_r \psi_j(t)^2 \mathrm{~d} t<\infty$ for some $r \in \mathbb{N}_+$.
\end{assumption}
\vspace*{5pt}
\vspace*{5pt}
\begin{assumption} \label{Ath} $\theta \rightarrow 0$  as  $n \rightarrow \infty$.
\end{assumption}
\vspace*{5pt}
\noindent Under the above assumptions, we now present the main result of the section.
\vspace*{10pt}
\begin{proposition} \label{Pcons}
{\rev Assume further that $\|\hat{\mathsf{K}}_{\mathsf{S}^2(\mathbb{X}^{[q]})}-\mathsf{K}_{\mathbb{X}^{[q]}}\|_{\mathrm{HS}}\to0$ a.s.} For each \bl{$1\leqslant j\leqslant q$}, if $\kappa_j$ has multiplicity $\mathcal{y}_j$,
\begin{enumerate}
\item there exists $\mathcal{y}_j$ sequences $(\hat{\kappa}_{\theta,l})_{l \in \mathcal{Y}_j}$ converging to $\kappa_j$, a.s.
\item $\mathsf{P}_{\hat{\psi}_{\theta,j}}$ converges to $\mathsf{P}_{\psi_j}$ in \bl{$\mathcal{B}_{2,H^{[q]}}$}, a.s.
\item \bl{if $\mathcal{y}_j=1$ then $\min_{s\in\{-1,1\}}\bigl\|s\,\hat{\psi}_{\theta, j}/\|\hat{\psi}_{\theta, j}\|-\psi_j\bigr\|\to0$ in $H=L^2_{[0,1]}$,  a.s.}
\end{enumerate}
\end{proposition}
\noindent {\rev Given this convergence, the conclusions follow from the perturbation theory of isolated eigenvalues of finite multiplicity \citep[Ch.~VIII, \S3]{kato80}. For fixed $q$ the convergence itself holds under Assumption \ref{Ath}, $\mathbb{E}\|\mathbb{X}^{[q]}\|^4<\infty$ and a.s. consistency of the estimated whitening operator, since $\|\mathsf{S}-I_{H^{[q]}}\|=\mathcal{O}(\theta)$ on $H^{[q]}$.} See also \cite{Silverman96,Boente00,Permantha17}, which provide similar results for the smoothed FPCA. A corollary of the proposition is that the subspace associated with the lowest eigenvalue orthogonal to the eigenspaces associated with normalized eigenvalues $\geq 3$ of multiplicity $\mathcal{y}_j$, is a consistent estimator for Fisher’s subspace\bl{, provided $\mu^{[q]}_\Delta\neq0$ and $\pi_0(1-\pi_0)>1/6$}. \cite{Li16} proposed an alternative type of consistency requiring distinct kurtoses eigenvalues to guarantee the independence of the corresponding projections. In our approach, this requirement is relaxed, as the focus is on estimating Gaussian components.

\section{Estimation with basis expansions}
\label{Compu.det}
\noindent
A general strategy for solving the continuous eigenproblem to an equivalent matrix eigenanalysis is to consider a representation of the empirical counterpart of $X$ with a finite basis of functions. Let $X^{[q]}(t)=(X_1^{[q]}(t), \ldots, X_n^{[q]}(t))^{\top}$ be a vector-valued function containing $n$ copies of $X$ assumed in a $q$-dimensional Hilbert space. Each function of $X^{[q]}(t)$ admits the basis function representation
\begin{equation}
X^{[q]}(t)=A\phi(t),
\label{expan}
\end{equation}
 where  $A\in \mathbb{R}^{n\times q}$ is a matrix of coefficients and $\phi(t)=(\phi_1(t),\ldots,\phi_q(t))^\top$ their respective vector of basis functions. The linear span of $\phi(t)$ is denoted by $H^{[q]}$ with inner product defined as $\left\langle f,g\right\rangle =\mathtt{f}^{\top}\mathcal{G}\mathtt{g}$, where $\mathtt{f},\mathtt{g}$ are the coefficient vectors of the functions $f,g\in H^{[q]}$ and  $\mathcal{G}=\langle \phi_{j},\phi_{k}\rangle \in \mathbb{R}^{q\times q}, j,k\in \{1\ldots q\},$ that is, the inner products of each pair of basis functions, so that possibly $\mathcal{G}\neq I_{q}$ when $\phi(t)$ may not be orthonormal. The first step in functional ICA is to pre-whiten the data, thus assume  $\mathbb{X}^{[q]}(t)=\tilde{A}\phi(t)$ be a set of whitened functional data, i.e., a basis expansion \bl{whose coefficient matrix satisfies $n^{-1}\mathcal{G}^{1/2}\tilde{A}^\top\tilde{A}\mathcal{G}^{1/2}=I_q$, equivalently $n^{-1}\tilde{A}^\top\tilde{A}=\mathcal{G}^{-1}$}. Computational algorithms for whitening functional data are provided in \cite{Vidal22}. Then, from expression (\ref{eq:stkurt}), we can define the sample kurtosis operator of  $\mathbb{X}^{[q]}$ as 
\begin{equation*}
\begin{split}
    \hat{\mathsf{K}}_{\mathbb{X}^{[q]}}(f)(s)	&=n^{-1}\sum_{i=1}^{n}\left\langle \mathbb{X}_{i}^{[q]},\mathbb{X}_{i}^{[q]}\right\rangle \left\langle \mathbb{X}_{i}^{[q]},f\right\rangle \mathbb{X}_{i}^{[q]}(s)=\left\langle n^{-1}\sum_{i=1}^{n}\left\Vert \mathbb{X}_i^{[q]}\right\Vert ^{2}\mathbb{X}_{i}^{[q]}(s)\mathbb{X}_{i}^{[q]},f\right\rangle \\
	& =\langle \hat{K}^{[q]}(s,\cdot),f\rangle,
\end{split}
\end{equation*}
where $\hat{K}^{[q]}(s,t)$ is a kurtosis kernel function admitting the following representation in terms of an orthonormalized basis 
\begin{equation*}
\hat{K}^{[q]}(s,t)=\phi^\top(s)\mathcal{G}^{-1/2}(n^{-1}\mathcal{G}^{1/2}\tilde{A}^\top D\tilde{A}\mathcal{G}^{1/2})\mathcal{G}^{-1/2}\phi(t),
\end{equation*}
{\rev where $D=D_0=\mathrm{diag}(\tilde{A}\mathcal{G}\tilde{A}^\top)$, i.e. $D_{0,ii}=\|\mathbb{X}_{i}\|^{2}$. The penalized problem of Proposition \ref{Prop4} instead carries the weights
\begin{equation*}
D_\theta=\mathrm{diag}(\tilde{A}\tilde{\mathcal{G}}_\theta\tilde{A}^\top), \qquad
\tilde{\mathcal{G}}_\theta=\mathcal{G}\mathcal{G}_\theta^{-1}\mathcal{G}=(L^{-1}\mathcal{G})^\top(L^{-1}\mathcal{G}), \qquad
\mathcal{G}_\theta=\mathcal{G}+\theta P^r=LL^\top,
\end{equation*}
i.e. $D_{\theta,ii}=\|\mathsf{S}(\mathbb{X}_{i})\|^{2}$, which reduces to $D_0$ at $\theta=0$.}

\vspace*{10pt}
\begin{proposition} \label{Prop3}
	Given the basis expansion in \eqref{expan}, the functional ICA of $X_i^{[q]}$ with respect to the  inner product $\langle \cdot,\cdot\rangle$ is equivalent to the  multivariate ICA of  matrix $A\mathcal{G}^{1/2}$ with the usual  metric in $\mathbb{R}^q.$
\end{proposition}

\vspace*{10pt}
\begin{proposition} \label{Prop4}
For any $\theta>0$, the penalized functional ICA of $X_i^{[q]}$ defined by the successive optimization problem in \eqref{1PSV} is equivalent to the multivariate ICA of the matrix $A\mathcal{G}^{1/2}$ using the metric $\mathcal{M}= (L^{-1} \mathcal{G}^{1/2})^\top (L^{-1} \mathcal{G}^{1/2})$ in $\mathbb{R}^q,$ with $L$ defined by the factorization $\mathcal{G}_\theta=\mathcal{G}+\theta P^r = LL^\top,$ and $P^r$ the matrix  whose elements are $\langle \mathsf{D}_r\phi_j,\mathsf{D}_r\phi_k\rangle.$
\end{proposition}

\vspace*{10pt}
From the relation between inner products given by \eqref{rel}, it can be deduced that the operator $\mathsf{S}^2$ is defined as $\mathsf{S}^{2}(f)=\phi(t)^{\top}\left(\mathcal{G}+\theta P^r\right)^{-1} \mathcal{G} \mathtt{f}$, with $f=\phi(t)^\top \mathtt{f}$. Then, for the smoothed whitened data $\mathsf{S}(\mathbb{X})$, the independent component scores are obtained as \bl{$\hat{\xi}_{\theta,j}=\tilde{A}\mathcal{G}(L^{-1})^\top v_{\theta,j}$}, and the kurtosis eigenfunctions as $ \hat{e}_j=\mathsf{S}^{-1}(\hat{\psi}_{\theta,j}).$ The smoothed IC scores do not have identity covariance for $\theta > 0$, but a transformation can be applied following Proposition 2 in \cite{Vidal21}.

\section{Numerical experiments}
\subsection{Simulated data}
To investigate the empirical performance of the proposed estimators, we conduct a study that extends the results of Simulation 2 in \cite{Pena14} using three possible taxonomies of mean differences in binary classification of functional data. Let  $X_i \ (i=1,\dots,n)$ be a mixture of two subpopulations  $\Pi_k \ (k=0,1),$  with  $n_k=n/2$ curves sampled on a grid of 20 equispaced points on $t \in [1,T]$ with $T=20$. Both groups have same quadratic covariance matrix cov$(t_j,t_{k})=\exp\{ -(2\ell^{2})^{-1}(t_{j}-t_{k})^{2}\}, \ j,k\in \{1,\dots,T\},$ with $\ell=15$. The data is then generated as
\begin{equation*}
X_{ij}=\sum_{k=0}^1\left(\sum_{j=1}^T \lambda_{j}^{1 / 2} Z_{k,ij} \gamma_{j}+\mu_{k} +\epsilon_{k,ij}\right) \mathbb{I}\left(X_i \in \Pi_k\right),
\end{equation*}
where $Z_{i j}$ are Gaussian random variables, $\epsilon_{k,ij}$ is an additive error term and $\mathbb{I}$ denotes the indicator function. Further extensions of the above model to non-Gaussian settings using  $Z_{k,i j} \sim \operatorname{exp}(1)-1$ can be also found in our results.

We consider the following versions of the above model: in Example 1 we define $\mu_0=0$ and $\mu_1=0.2 \cos (3 \pi t/T)$, the means differ in shape; in Example 2, $\mu_0=0.3 \cos (3 \pi t/T)$ and $\mu_1=0.2 \cos (3 \pi t/T)$, the means have equal shape and slightly differ in amplitude; in Example 3 we set $\mu_0=0.2 \sin (3 \pi t/T)$ and $\mu_1=0.2 \cos (3 \pi t/T)$, the means are equal in shape but dephased $\pi/2$. In all cases, $Z_{k,i j}$ are sampled from a standard normal distribution and $\epsilon_{k,ij}\sim N(0, \sigma^2)$.  We generated 200 datasets for each experiment with sample sizes $n_k=30,50$. The choice of $\beta$ includes both functional PCA estimators, such as the dominant eigenfunction of the covariance operator and the one that minimizes the kurtosis of the projection scores, as well as functional ICA estimators, which are based on various whitening procedures and regularization schemes (for which we use the last kurtosis eigenfunction). As in many other studies, the roughness penalty is defined as the integrated squared derivative of order $r=2$. The \textsf{R} package \textsf{pfica} \citep{pfica} was used for the implementation of various functional pre-whitening methods via B-spline expansions with $q=5$.

Results for $\sigma=0$ are shown in Table~\ref{tab1}. In all examples, the overall good behaviour of kurtosis classifiers based on the last independent component (minimum kurtosis) is apparent, particularly for the smoothed kurtosis projections and large sample sizes. In Example 3, the PC with the lowest kurtosis coefficient performed notably well, similarly to their kurtosis peers and eventually outperforming the rates of the non-smoothed kurtosis. Regarding functional whitening, results indicate that classification optimization with the proposed operators is not that different, although Cholesky  whitening reaches good performance in Examples  1 and 3, while zero-phase components analysis whitening does better in Example 2. In the non-Gaussian simulation, results are more balanced between both functional ICAs, although superior to the rest of competitors.
 
\begin{figure}[h]
\centering\includegraphics[scale=.30]{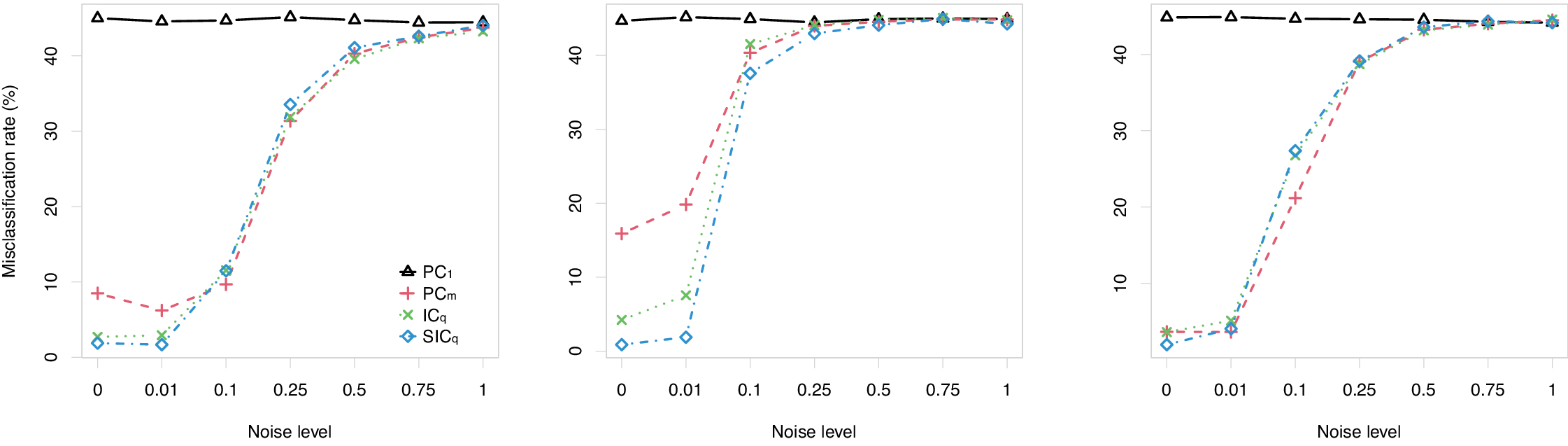}
\caption{Comparison of the misclassification rate \bl{(mean values over the 200 replications)} for the three Examples 1-3 (from left to right) and different levels of noise ($\sigma$) using ZCA whitening. PC$_1$, first principal component; PC$_m$, principal component with lowest kurtosis coefficient; IC$_q$, $q$th independent component (minimal kurtosis); SIC$_q$, $q$th smoothed independent component.}
\label{fig0}
\end{figure}

Fig.~\ref{fig0} further illustrates the effect of modulating the noise in the performance of the classifiers on a Gaussian scenario. Note that as $\sigma$ grows, it exponentially worsens the classification rate. 

Notwithstanding, results for the smoothed kurtosis are very competitive for mild levels of noise. These analyses point to the importance of finding a good trade-off between groups when smoothing the data, as both noise and the type of smoothing (especially if it is homogeneous across curves) can undermine the  effectiveness of these classifiers.

\begin{table}[h!]
\caption{ Simulation results for the mean and standard deviation (in parentheses) of the classification errors obtained with 200 repetitions of the experiment for different sample sizes and zero error variance. PC$_1$, first principal component; PC$_m$, principal component with lowest kurtosis coefficient; IC$_q$, $q$th independent component (minimal kurtosis); SIC$_q$, $q$th smoothed independent component; PCA, principal component analysis whitening; PCA-cor, principal component analysis correlated whitening; ZCA, zero-phase component analysis or Mahalanobis whitening; ZCA-cor, zero-phase component analysis or Mahalanobis correlated whitening; Cholesky, Cholesky whitening.}
\label{tab1}
\vskip-0.3cm\hrule

\smallskip
\centering\tiny

\begin{tabular}{ccccccc}

     &       & \multicolumn{5}{l}{Results (\%) for the centroid classifiers:} \\ 
     \cmidrule{2-7}
Data & $n_k$ & PC$_1$ & PC$_m$ & Whitening & IC$_q$ & SIC$_q$  \\ 
      \midrule
       &       & Scenario I (Gaussian) \\ 
Example 1 & 30 & 45.49 (3.295) & 6.167 (10.43) & PCA & 2.583 (3.698) & 2.520 (5.193) \\ 
   &   &   &   & PCA-cor & 2.583 (3.698) & 2.540 (5.869) \\ 
   &   &   &   & ZCA & 2.592 (3.709) & 1.590 (2.632) \\ 
   &   &   &   & ZCA-cor & 2.592 (3.709) & 1.955 (3.234) \\ 
   &   &   &   & Cholesky & 2.575 (3.702) & 1.219 (2.083) \\ 
   & 50 & 45.47 (3.382) & 3.985 (7.754) & PCA & 0.970 (1.662) & 0.588 (1.058) \\ 
   &   &   &   & PCA-cor & 0.970 (1.662) & 0.885 (2.122) \\ 
   &   &   &   & ZCA & 0.965 (1.661) & 0.531 (1.051) \\ 
   &   &   &   & ZCA-cor & 0.970 (1.671) & 0.678 (1.066) \\ 
   &   &   &   & Cholesky & 0.965 (1.661) & 0.490 (1.266) \\ 
  Example 2 & 30 & 44.83 (4.103) & 17.35 (15.10) & PCA & 4.500 (4.971) & 1.716 (2.318) \\ 
   &   &   &   & PCA-cor & 4.492 (4.960) & 1.766 (2.411) \\ 
   &   &   &   & ZCA & 4.475 (4.969) & 0.876 (1.431) \\ 
   &   &   &   & ZCA-cor & 4.475 (4.969) & 0.903 (1.455) \\ 
   &   &   &   & Chol & 4.492 (4.966) & 1.688 (2.342) \\ 
   & 50 & 45.98 (3.047) & 11.35 (11.32) & PCA & 1.750 (2.198) & 0.758 (1.308) \\ 
   &   &   &   & PCA-cor & 1.755 (2.200) & 0.969 (1.451) \\ 
   &   &   &   & ZCA & 1.755 (2.200) & 0.397 (0.790) \\ 
   &   &   &   & ZCA-cor & 1.750 (2.203) & 0.420 (0.787) \\ 
   &   &   &   & Cholesky & 1.740 (2.202) & 0.942 (1.517) \\ 
  Example 3 & 30 & 45.16 (3.810) & 2.975 (4.474) & PCA & 3.525 (4.031) & 1.667 (2.472) \\ 
  &   &   &   & PCA-cor & 3.533 (4.026) & 1.988 (4.355) \\ 
   &   &   &   & ZCA & 3.533 (4.026) & 2.264 (4.476) \\ 
   &   &   &   & ZCA-cor& 3.525 (4.018) & 2.004 (4.214) \\ 
   &   &   &   & Cholesky & 3.517 (4.023) & 1.286 (2.117) \\ 
   & 50 & 45.54 (3.230) & 1.910 (1.560) & PCA & 1.530 (1.974) & 0.941 (1.541) \\ 
   &   &   &   & PCA-cor & 1.535 (1.974) & 1.129 (1.775) \\ 
   &   &   &   & ZCA & 1.535 (1.974) & 1.043 (1.700) \\ 
   &   &   &   & ZCA-cor & 1.535 (1.974) & 0.867 (1.505) \\ 
   &   &   &   & Cholesky & 1.530 (1.974) & 0.663 (1.044) \\
    &       &  \\ 
    \midrule
     &       & Scenario II (non-Gaussian) \\ 

Example 1 & 30 & 45.40 (3.592) & 6.742 (10.45) & PCA & 3.917 (3.967) & 3.909 (5.068) \\ 
   &   &   &   & PCA-cor & 3.917 (3.967) & 4.167 (5.950) \\ 
   &   &   &   & ZCA & 3.900 (3.970) & 3.438 (2.643) \\ 
   &   &   &   & ZCA-cor & 3.900 (3.970) & 3.452 (2.602) \\ 
   &   &   &   & Cholesky & 3.908 (3.967) & 3.563 (2.543) \\ 
   & 50 & 45.97 (3.271) & 4.875 (9.073) & PCA & 3.075 (3.547) & 2.920 (3.190) \\ 
   &   &   &   & PCA-cor & 3.075 (3.547) & 2.819 (2.889) \\ 
   &   &   &   & ZCA & 3.075 (3.547) & 3.010 (3.448) \\ 
   &   &   &   & ZCA-cor & 3.075 (3.547) & 3.071 (3.936) \\ 
   &   &   &   & Cholesky & 3.075 (3.547) & 2.678 (2.864) \\ 
  Example 2 & 30 & 45.208 (4.058) & 7.45 (10.58) & PCA & 8.250 (9.897) & 8.349 (10.48) \\ 
   &   &   &   & PCA-cor & 8.250 (9.897) & 8.063 (10.77) \\ 
   &   &   &   & ZCA & 8.250 (9.889) & 7.155 (9.419) \\ 
   &   &   &   & ZCA-cor & 8.242 (9.892) & 6.800 (9.561) \\ 
   &   &   &   & Cholesky & 8.233 (9.897) & 7.343 (9.434) \\ 
   & 50 & 46.29 (2.748) & 7.805 (11.97) & PCA & 7.145 (10.48) & 6.053 (8.575) \\ 
   &   &   &   & PCA-cor & 7.150 (10.50) & 5.840 (8.591) \\ 
   &   &   &   & ZCA & 7.155 (10.50) & 6.694 (10.56) \\ 
   &   &   &   & ZCA-cor & 7.155 (10.50) & 6.640 (10.69) \\ 
   &   &   &   & Cholesky & 7.150 (10.48) & 6.075 (8.519) \\ 
  Example 3 & 30 & 45.15 (4.049) & 10.92 (11.87) & PCA & 6.858 (7.936) & 5.951 (6.266) \\ 
   &   &   &   & PCA-cor & 6.858 (7.936) & 6.516 (7.232) \\ 
   &   &   &   & ZCA & 6.867 (7.934) & 7.060 (8.419) \\ 
   &   &   &   & ZCA-cor & 6.867 (7.934) & 6.425 (7.626) \\ 
   &   &   &   & Cholesky & 6.858 (7.936) & 5.654 (6.642) \\ 
   & 50 & 45.94 (3.092) & 10.92 (11.69) & PCA & 5.085 (6.098) & 4.832 (5.616) \\ 
   &   &   &   & PCA-cor & 5.085 (6.098) & 5.119 (6.249) \\ 
   &   &   &   & ZCA & 5.095 (6.120) & 4.935 (6.120) \\ 
   &   &   &   & ZCA-cor & 5.100 (6.119) & 4.873 (5.696) \\ 
   &   &   &   & Cholesky & 5.095 (6.095) & 4.173 (4.228) \\ 
\midrule

\end{tabular}
\end{table}

\subsection{Real datasets} \label{ss5.2}
\begin{figure}[h]
\centering\includegraphics[scale=.28]{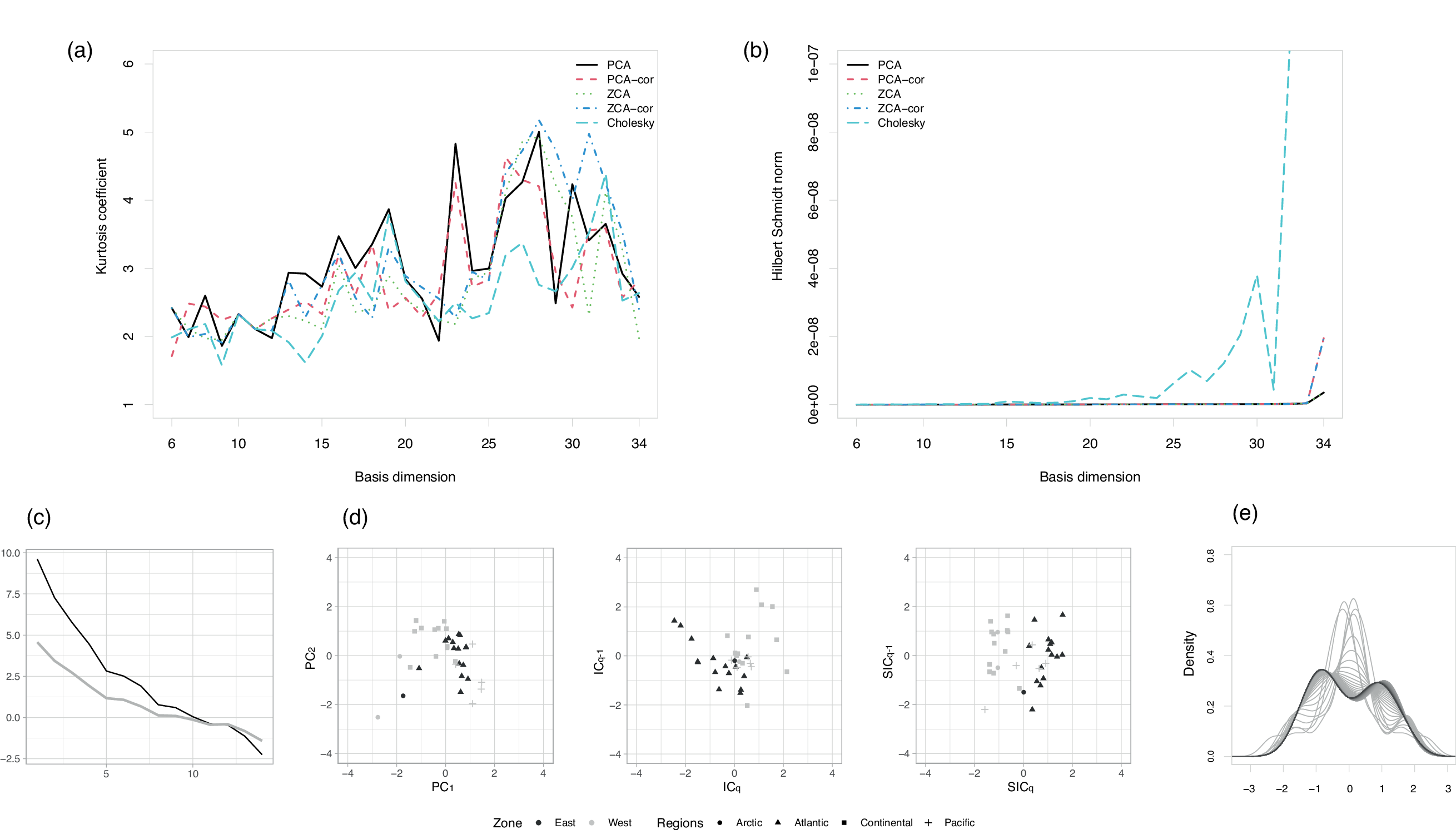}
\caption{ Canadian Weather data results. (a) Kurtosis coefficient of  $\hat{\xi}_{\theta,q} \ (\theta=0,100,\dots,3\cdot10^4),$ for  various B-spline basis dimensions ($q$) and whitening procedures. (b) The plot shows the asymptotic behavior of the whitening transformations when $q$ grows up to $n$ as evaluated by \bl{$\|\Delta \hat{\Gamma}_{\mathbb{X}^{[q]}}\|_{\mathrm{HS}}$} \ $(q=6,\dots,34)$ (c) Picard's plot. The black line stands for the eigenvalues $\log(\hat{\lambda}_j)$ and the grey one, for the means of absolute values of principal component scores given by $\sum_{i=1}^n |\langle X^{[q]}_i,\hat{\gamma}_j \rangle|$ expressed in a logarithmic scale. (d) Scatter plots. From left to right: functional PCA, ICA and smoothed functional ICA using a basis expansion of $q=14$ (e) Estimated densities of the vector $\hat{\xi}_{\theta,q}$ for each lambda, showing the effect of smoothing the kurtosis operator.}
\label{Fig1}
\end{figure}

Our methods are now applied to well-known datasets in the functional data literature.  In the first example, we show that the smoothed kurtosis is able to find bimodality in the Canadian Weather data, which is usually treated as a discrimination problem of more than two groups. We consider a geographical division based on a west-east location distribution rather than the usual four climate regions. The whitening method and penalty parameter was selected using cross-validation by minimizing the kurtosis coefficient of the projections on to $\hat{\psi}_q$ (Fig. \ref{Fig1}a) with $q\in\{5,\dots,34\}$, $\theta\in\{0,100,\dots,10^8\}$. Results suggest the presence of bimodality in these data, and the few misclassified observations appear to be locations close to large bodies of water, commonly encountered in the west zone. \bl{Due to the} representativeness of these data, we further asses the asymptotic behavior of \bl{$\|\Delta\hat{\Gamma}_{\mathbb{X}^{[q]}}\|_{\mathrm{HS}}=\| \hat{\Gamma}_{\mathbb{X}^{[q]}}- I_{H^{[q]}}\|_{\mathrm{HS}}$} as $q$ increases, to verify the consistency of the proposed methods.  As shown in Fig. \ref{Fig1}b, most of the whitening procedures  (under mild conditions, i.e. $q<n$) converge to the true parameter, with Cholesky whitening being the less consistent approach, albeit the one that provides more interesting results. 

\begin{figure}[h]
\centering\includegraphics[scale=.28]{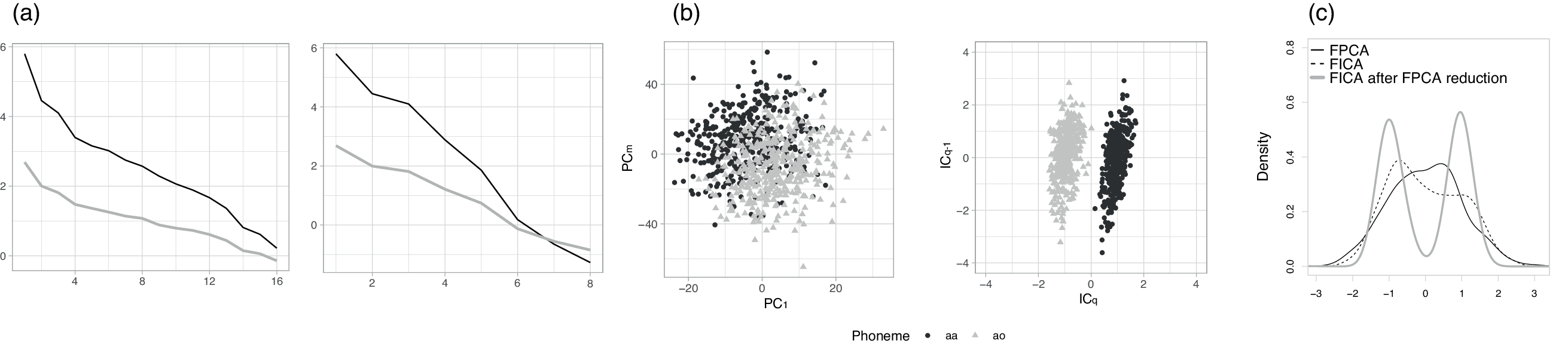}
\caption{Phoneme data results. (a) Picard's plot before and after a functional PCA per group. (b) Scatter plots. From left to right: functional PCA (first component against component with lowest kurtosis), ICA (non-smoothed) after reduction. (c) Estimated densities of the vector of scores with lowest kurtosis using different reduction techniques.}
\label{Fig2}
\end{figure}

In a second example, we consider the phoneme dataset as analyzed in \cite{Delaigle12}. The data were retrieved from the \textsf{fds} package \citep{fds} and consist of 400 log-periodograms constructed from audio recordings of males pronouncing the phonemes `aa' as in dark and `ao' as in water. The similarity between both groups of curves has been previously reported to pose a challenging problem of classification. In fact, we were neither able to find interesting projections with any of the proposed methods. As workaround, we propose to perform a functional PCA on each sample and use the basis function expansion
$
X^{[p]}_{k,i}= \sum_{j=1}^p \langle X^{[q]}_{k,i},\hat{\gamma}_{k,j} \rangle \hat{\gamma}_{k,j}\mathbb{I}\left(X_i \in \Pi_k\right),
$
which takes the matrix form
\begin{equation*}
X^{[p]}_{k} = (A_{k}\mathcal{G}b_{k}^\top) b_{k}\phi (t),
\end{equation*}
where $b_{k}=U_{k}\mathcal{G}^{-1/2}$, with $U_{k} \in \mathbb{R}^{p\times q}$ the matrix of eigenvectors of $n^{-1}\mathcal{G}^{1/2}A^\top_{k} A_{k} \mathcal{G}^{1/2}$ truncated at the $p$-row and $X^{[p]}_{k}=(X^{[p]}_{k,1},\dots,X^{[p]}_{k,n_k})^\top$. Taking $X^{[p]}=(X^{[p]}_{1},X^{[p]}_{2})$  and using the coefficients in terms of basis functions pooled by rows, one can perform the functional ICA using these coefficients. The matrix $A=\{(A_{k}\mathcal{G}b_{k}^\top) b_{k}\}_{k=0}^1$, however, might have non invertible covariance matrix. Although this could  be reversed by a suitable Tikhonov regularization, to avoid harming the whitening procedure, the best option is to truncate and perform the functional ICA on the $p$-principal components. Therefore, we use the first components up to the limit where the whitening transformation no longer meets the orthonormality property. Results in Fig. \ref{Fig2}, taking $p=8$ components, show the great improvement of performing functional ICA on these representations, which achieves near-perfect classification with an error rate of  0.125 \%. Applying smoothed functional ICA after the functional FPCA reduction did not improve the results.

\section{Modeling spectral biomarkers for depressive disorder}
As we have seen, having a large sample size can be advantageous for classifying data, although its success depends on several other factors. In this section we focus on \bl{relatively} small sample sizes with hypervariability and possibly heterogeneous variances. This problem is recurrently found in neuroscientific data and in situations where a consistent data size from the treatment group is not available. Here, we use the publicly accessible Leipzig Mind-Brain-Body (LEMON) dataset for modeling electroencephalographic (EEG) biomarkers for depressive disorder. Prior studies investigating these data have focused on modeling neuronal signatures typically associated with \bl{Parkinson's} disease \citep{Zhang21}, acute psychological stress \citep{Reinelt19} and aging in the human brain \citep{Kumral20}. Thorough descriptions of the data collected and experimental procedure can be found in \cite{Babayan19}. Briefly, 221 healthy participants, young and older adults, were recruited to investigate mind-body-emotion interactions across a wide range of cognitive, emotional and physiological phenotypes. Among all psychological and physiological assessment conducted, we are interested in the Hamilton depression scale (HDS) scores and respective participant's EEG measurements. From the total population, we consider the group that scored higher (from 8 to 14) on HDS (11 participants, 7 females), to compare them to a subpopulation of equal size randomly selected from participants scoring 0 (49 participants, 30 females - control group).

The analyzed data consists of 62-channel EEG of 8 minutes resting-state  recordings with eyes closed. Signals were pre-processed using common techniques for reducing artifactual activity, resampled at 250 Hz and bandpass filtered within 1-45 Hz using an 8th-order Butterworth filter \citep{Babayan19}. As biomarker, we use a broad-band Log power spectral density (5-45 Hz) of the EEG time series calculated using the Welch transform (Hamming window of 20 samples, 50\% overlap). The curves (grouped per channel) were approximated using B-spline expansions $(q=11)$ for subsequent analyses. In all experiments, we consider equal sample sizes for both groups, even in the case of missing channels. For a fixed treatment group, 1000 iterations bootstrap were performed by randomly selecting a subpopulation of the control group without repetition. Results were averaged. We asses the performance of functional PCA (minimum kurtosis), smoothed functional ICA (ZCA whitening, $\theta \in \{0,1,\dots 10\}$), also with the functional representations proposed in Sect.~\ref{ss5.2}. For each channel, the smoothing parameter was selected as the one that minimized the misclassification rate via cross-validation on the last ICs. After reduction, only 20 out of 62 cases showed improvement with the smoothed estimators. The tolerance threshold of the covariance eigenvalues when truncating at $p$ was set to $5\times 10^{-10}$.

\begin{figure}[h!] 
\centering\includegraphics[scale=.25]{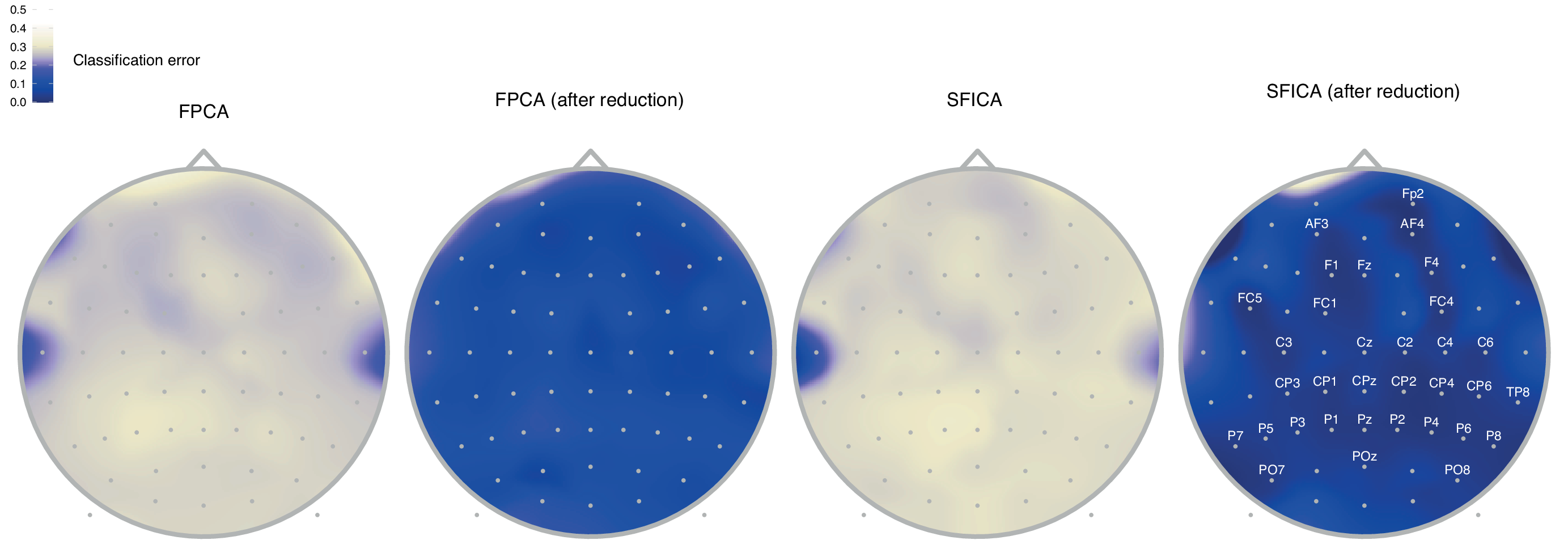}
\caption{Topographic scalp heatmaps display the classification error for the assessment of the log power spectral density curves at each sensor location, with each channel represented by grey points. The classification error scores are interpolated across the circular field to improve interpretability. Channels scoring a classification error $< 0.05$ are labeled.}
\label{Fig3}
\end{figure}

Fig. \ref{Fig3} shows the classification error for each channel through interpolated topographical maps. Results clearly indicate that the functional representations proposed in Sect.~\ref{ss5.2} perform favorably compared to when data is not smoothed per group. More concretely, the functional smoothed ICA outperforms the functional PCA after the reduction.  Consistent with previous research, channel locations that allow to optimally classify new observations in the medial frontal, Fp2 area and sensory motor cortices have been previously associated with depression pathology; see, for example, \cite{Bludau16,Cook14,Ray21}. The temporal lobe region (the lateral blue zones of the FPCA/SFICA results) has also been reported to exhibit physiological changes in patients with depression disorder due to a thinning of the cortical thickness \citep{Schmaal16}. Accordingly, our results suggest that a decrease of the misclassification rate in these areas  could possibly be  related to spurious muscular activity rather than genuine brain activity, as these effects seem to disappear (or to attenuate) after the reduction.

\section*{Concluding remarks}
In this paper, we developed a functional ICA method based on penalizing the roughness of the kurtosis eigenfunctions by its integrated squared derivative as a method for regularization. We have shown this method has competitive operating features in binary classification problems, both in Gaussian and non-Gaussian settings. While the proposed form of regularization positively impacts the classification performance and consistency, the revised real dataset examples show that its benefits do not always surpass the effectiveness of a marginal KL reduction in cases of complex variability. Regularization is useful when the conditions for differentiability are met such that the last (normalized) kurtosis eigenvalue converges toward 1, allowing a targeted selection of the penalty parameter through cross-validation. Careful consideration should be given to the deliberate application of cross-validation without ensuring adherence to the consistency criteria established in Sect.~\ref{SCons}, as this could lead to unreliable estimates. It would be interesting to study how, in the case of a fractional differential operator, adjusting this parameter might improve results or offer deeper insights into consistency.

It has been further shown that a subset of the eigenfunctions of the kurtosis operator corresponds to Fisher’s discriminant subspace under certain Gaussian assumptions, and we have established the consistency of these eigenfunctions as estimators for the discriminant directions in the population setting. However, as pointed out in \cite{Delaigle12} ``the general case seems not to admit an elementary, insightful derivation of the optimal $\beta$” which rather seems to indicate, based on prior research \citep{Baillo10,Berrendero18}, that the Feldman-H\'{a}jek dichotomy might be the motor underlying functional classification problems. \bl{Building on this result}, our approach recasts on ``Gaussianizing”  the data using a suitable whitening method, regularization or a group-wise functional PCA reduction in order to balance and minimize the kurtosis of the eigenprojections.  In conclusion, the kurtosis operator (as a measure of non-Gaussianity), and particularly its spectral attributes, \bl{offers an analytical pathway that can bring us closer to achieving near-perfect accuracy}, as elucidated by the Feldman-H\'{a}jek dichotomy.

\appendix
\titleformat{\section}[hang]{\large\normalfont\scshape\filcenter}{Appendix~\thesection.}{1em}{}

\section{Proofs of formal statements} \label{secA0}
\begin{proof}[Proof of Proposition 2.1]
For $\theta \geq 0$, consider the operator $\mathsf{S}$ has spectral decomposition
\begin{equation*}
\bl{\mathsf{S}=(I_{H}+\theta  \mathsf{D}_r^*\mathsf{D}_r)^{-1 / 2}}=\bl{\int\left(1+\theta\mathfrak{d}\right)^{-1 / 2} \mathrm{d} E(\mathfrak{d})},
\end{equation*}
{\rev
\noindent where $E$ is the resolution of the identity of $\mathsf{D}_r^*\mathsf{D}_r$. Since $\mathfrak{d}\geq0$ and $\theta \geq 0$, it follows that $0<(1+\theta\mathfrak{d})^{-1/2} \leq 1$, so that the operator norm $\left\|\mathsf{S}\right\|_{\mathcal{B}_H}$, which is the supremum of the absolute values of these eigenvalues, satisfies $\left\|\mathsf{S}\right\|_{\mathcal{B}_H} \leq 1$. Boundedness alone does not give the trace-class property. Since $\operatorname{tr}(g\otimes g)=\|g\|^2$ for $g \in H$, we have
\begin{equation*}
\operatorname{tr}\mathsf{K}_{\mathsf{S}(\mathbb{X})}=\mathbb{E}\|\mathsf{S}(\mathbb{X})\|^4,
\end{equation*}
so that the statement is equivalent to Assumption \ref{finite4}. In the Gaussian case it holds for every $\theta>0$ such that $\operatorname{tr}(\mathsf{S}^2)=\sum_{j}(1+\theta\mathfrak{d}_j)^{-1}<\infty$, where $(\mathfrak{d}_j)$ are the eigenvalues of $\mathsf{D}_r^*\mathsf{D}_r$, since then $\mathbb{E}\|\mathsf{S}(\mathbb{X})\|^4=\{\operatorname{tr}(\mathsf{S}^2)\}^2+2\operatorname{tr}(\mathsf{S}^4)<\infty$.
}
See  also \cite{Permantha17} for further details on the operator $\mathsf{S}$. 

\end{proof}

\begin{proof}[Proof of Proposition 2.2]
{\rev
By \eqref{rel}, $\langle f,\mathsf{S}^2(\mathbb{X})\rangle_\theta=\langle f,\mathbb{X}\rangle$ and $\langle \mathsf{S}^2(\mathbb{X}),\mathsf{S}^2(\mathbb{X})\rangle_\theta=\langle \mathbb{X},\mathsf{S}^2(\mathbb{X})\rangle=\|\mathsf{S}(\mathbb{X})\|^2$, whence
\begin{equation*}
\mathsf{K}_{\mathsf{S}^2(\mathbb{X})}(f)=\mathbb{E}\left\{\langle f,\mathbb{X}\rangle\|\mathsf{S}(\mathbb{X})\|^2\mathsf{S}^2(\mathbb{X})\right\}=\mathsf{S}^2\mathsf{K}_\theta(f).
\end{equation*}
Since $\langle f,\mathsf{K}_\theta\psi\rangle=\langle f,\mathsf{S}^2\mathsf{K}_\theta\psi\rangle_\theta$, the equation $\langle f,\mathsf{K}_\theta\psi_{\theta,j}\rangle=\kappa_{\theta,j}\langle \psi_{\theta,j},f\rangle_\theta$ for all $f$ is equivalent to $\mathsf{K}_{\mathsf{S}^2(\mathbb{X})}\psi_{\theta,j}=\kappa_{\theta,j}\psi_{\theta,j}$.
}
\end{proof}

\begin{proof}[Proof of Proposition 2.3]
Let us denote by $\mathsf{K}_{\mathsf{S}^2(\mathbb{X})}$ the kurtosis operator of the random variable $\mathsf{S}^2(\mathbb{X})$ with the inner product $\langle \cdot,\cdot\rangle_\theta$ and by $\mathsf{K}_{\mathsf{S}(\mathbb{X})}$ 
the kurtosis operator of $\mathsf{S}(\mathbb{X})$ with $\langle \cdot,\cdot\rangle.$

The equivalence between the two eigensystems is clearly deduced from the following relationship between both kurtosis operators: 
\begin{equation*}
\mathsf{K}_{\mathsf{S}^2(\mathbb{X})} = \mathsf{S} \mathsf{K}_{\mathsf{S}(\mathbb{X})}  \mathsf{S}^{-1}.
\end{equation*}
In fact,
\begin{align*}
\mathsf{K}_{\mathsf{S}^2(\mathbb{X})} (f) & =  \mathbb{E}\{\langle f,\mathsf{S}^2(\mathbb{X}) \rangle_\theta \langle \mathsf{S}^2(\mathbb{X}),\mathsf{S}^2(\mathbb{X}) \rangle_\theta \mathsf{S}^2 (\mathbb{X})\} \\
& =  \mathsf{S} [  \mathbb{E}\{\langle \mathsf{S}^{-1}(f),\mathsf{S}(\mathbb{X}) \rangle \langle \mathsf{S}(\mathbb{X}),\mathsf{S}(\mathbb{X}) \rangle \mathsf{S} (\mathbb{X})\} ] \\
&  =  \mathsf{S} \mathsf{K}_{\mathsf{S}(\mathbb{X})} ( \mathsf{S}^{-1} (f)).
\end{align*}
\end{proof}

{\rev
\begin{proof}[Proof of Corollary~\ref{cor2}]
We begin with the Fisher score before standardization. Let
$Z=\langle X-\mu,\beta^{[q]}\rangle$ with
$\beta^{[q]}=\sum_{j=1}^{q}\lambda_j^{-1}\mathtt{v}_j\gamma_j$. Since the two
classes have equal prior probabilities, $\mu=(\mu_0+\mu_1)/2$ and
$\mu_k-\mu=\pm\mu_\Delta/2$, so that
\begin{equation*}
\mathbb{E}(Z\mid Y=k)=\pm\tfrac12\langle\mu_\Delta,\beta^{[q]}\rangle=\pm\tfrac12 d_q^2,
\qquad
\operatorname{Var}(Z\mid Y=k)=\langle\beta^{[q]},\mathcal{R}\beta^{[q]}\rangle=d_q^2,
\end{equation*}
both equalities following from $\sum_{j\le q}\mathtt{v}_j^2/\lambda_j=d_q^2$, and the
variance being the same in the two classes because they share $\mathcal{R}$. By the
law of total variance,
\begin{equation*}
\operatorname{Var}(Z)=d_q^2+\frac{d_q^4}{4},
\qquad\text{so that}\qquad
\frac{d_q^4/4}{d_q^2+d_q^4/4}=\frac{d_q^2}{4+d_q^2}=p_q
\end{equation*}
is the share of $\operatorname{Var}(Z)$ due to the separation of the class means,
and $1-p_q$ is the share due to the common within-class variance.

Standardization leaves these two shares unchanged. Indeed, by Proposition \ref{Ke},
$\psi_q\propto\Gamma^{1/2\dagger}_{X^{[q]}}\mu^{[q]}_\Delta$, so
$\xi_q=\langle\mathbb{X}^{[q]},\psi_q\rangle$ is proportional to
$\langle X-\mu,\Gamma^{-1}_{X^{[q]}}\mu^{[q]}_\Delta\rangle$. Since
$\Gamma_{X^{[q]}}=\mathcal{R}_q+\tfrac14\mu^{[q]}_\Delta\otimes\mu^{[q]}_\Delta$, the
rank-one inverse formula gives
\begin{equation*}
\Gamma^{-1}_{X^{[q]}}\mu^{[q]}_\Delta
=\frac{1}{1+d_q^2/4}\,\mathcal{R}_q^{-1}\mu^{[q]}_\Delta
\propto\beta^{[q]},
\end{equation*}
so $\xi_q$ is proportional to $Z$ and standardization merely rescales $Z$ to unit
variance, without altering the two shares above.

The two contributions can now be identified explicitly. Since the classes are
Gaussian with equal priors, the class label contributes a symmetric two-point
component and the within-class fluctuation a centred Gaussian one, so
\begin{equation*}
\xi_q=\varepsilon\sqrt{p_q}+\sqrt{1-p_q}\,\zeta,
\end{equation*}
where $\varepsilon$ takes the values $-1$ and $1$, each with probability $1/2$,
indicating the class, and $\zeta\sim N(0,1)$ is independent of $\varepsilon$. In
particular $\mathbb{E}\xi_q^2=p_q+(1-p_q)=1$, as required by the standardization,
and
\begin{equation*}
\mathbb{E}\xi_q^4=p_q^2+6p_q(1-p_q)+3(1-p_q)^2=3-2p_q^2 .
\end{equation*}

We can now compute the spectrum of $\mathsf{K}_{\mathbb{X}^{[q]}}$. Since the total
covariance of $\mathbb{X}^{[q]}$ is $I_q$ and its between-class covariance is
$p_q(\psi_q\otimes\psi_q)$, its common within-class covariance is
$I_q-p_q(\psi_q\otimes\psi_q)$, so writing $\mathbb{X}^{[q]}=\xi_q\psi_q+W$ the
orthogonal component satisfies $W\sim N(0,I_{q-1})$ on $\psi_q^{\perp}$ and is
independent of $\xi_q$. The cross terms
$\mathbb{E}\{\|\mathbb{X}^{[q]}\|^2\xi_q\langle W,u\rangle\}$ vanish because $W$ is
symmetric, so $\mathsf{K}_{\mathbb{X}^{[q]}}$ is diagonal for the decomposition
$\operatorname{span}\{\psi_q\}\oplus\psi_q^{\perp}$, with
\begin{equation*}
\kappa_q=\mathbb{E}\xi_q^4+\mathbb{E}\|W\|^2\,\mathbb{E}\xi_q^2=\mathbb{E}\xi_q^4+(q-1),
\qquad
\kappa_u=\mathbb{E}\xi_q^2\,\mathbb{E}\langle W,u\rangle^2+\mathbb{E}\{\|W\|^2\langle W,u\rangle^2\}=1+(q+1)
\end{equation*}
for any unit $u\perp\psi_q$. Subtracting the Gaussian baseline $q-1$ gives
$\mathtt{kurt}_q=3-2p_q^2$ and $\kappa_u-(q-1)=3$.

It remains to connect these expressions with classification and with Gaussian measure
singularity. For the centroid classifier $T^0$ along
$\beta^{[q]}$ one has $\nu=\langle\mu_\Delta,\beta^{[q]}\rangle=d_q^2$ and
$\sigma_Q^2=\langle\beta^{[q]},\mathcal{R}\beta^{[q]}\rangle=d_q^2$, so
$\nu/(2\sigma_Q)=d_q/2$ and, $Q$ being Gaussian,
$\operatorname{err}_q=\Phi(-d_q/2)$. Finally $p_q=d_q^2/(4+d_q^2)$ is strictly
increasing in $d_q^2$, whence $\mathtt{kurt}_q\to1\iff p_q\to1\iff d_q^2\to\infty$;
and $d_q^2\to\infty$ if and only if
$\mu_\Delta\notin H(\mathcal{R})=\operatorname{ran}(\mathcal{R}^{1/2})$. Since $m_0$
and $m_1$ share the covariance operator $\mathcal{R}$, conditions 1 and 3 of
Theorem \ref{HF} hold automatically and the dichotomy is decided by this
Cameron--Martin condition, which gives \eqref{gchain}; the last equivalence
follows from $\Phi(-d_q/2)\to0\iff d_q\to\infty$.
\end{proof}
}

\begin{proof}[Proof of Proposition 4.1]
As indicated in Sect.~\ref{sec2}, the weight functions of the  independent components are obtained as the eigenfunctions of the sample kurtosis operator by solving the following eigenproblem:
\begin{equation*}
\hat{\mathsf{K}}_{\mathbb{X}^{[q]}} \psi(t)= \hat{\kappa} \hat{\psi}(t).
\end{equation*}
If we expand the independent component  weight functions as
$\hat{\psi}(t)=\phi(t)^\top b,$ then the problem turns in matrix form as
\begin{equation*}
n^{-1} \tilde{A}^{\top}\bl{D_0}\tilde{A} \mathcal{G} b = \hat{\kappa} b,
\end{equation*}
which is equivalent to
$$n^{-1} \mathcal{G}^{1/2}\tilde{A}^{\top}\bl{D_0}\tilde{A} \mathcal{G}^{1/2} u = \hat{\kappa} u,$$
with $u= \mathcal{G}^{1/2} b.$ That is, $\Sigma_{\tilde{A}\mathcal{G}^{1/2}}^{[4]} u = \hat{\kappa} u.$ Taking into account that the matrix $\tilde{A}\mathcal{G}^{1/2}$ is the whitening of the matrix $A\mathcal{G}^{1/2},$ we can conclude that the functional ICA is equivalent to ICA of matrix $A\mathcal{G}^{1/2}.$ 
\end{proof}

\begin{proof}[Proof of Proposition 4.2]
If we expand the  weight functions in (\ref{kurtpen}) as
$f(t)=\sum_{j=1}^{q}b_{j}\phi_{j}(t)=\phi(t)^\top b$
where $b=(b_1,\ldots,b_q)^\top,$ the coefficients of $\hat{\psi}_{\theta,j}$, are obtained by solving the penalized kurtosis problem (\ref{kurtpen}) expressed in matrix form as
\begin{equation*}
\frac{\bl{n^{-1}b^{\top}\mathcal{G}\tilde{A}^{\top} D_\theta\tilde{A}\mathcal{G}b}}{b^{\top}\mathcal{G}b+\theta b^{\top}P^rb}=\frac{b^{\top} \mathcal{G}^{1/2} \Sigma_{\tilde{A}\mathcal{G}^{1/2}}^{[4]} \mathcal{G}^{1/2} b}{b^{\top}\left(\mathcal{G}+\theta P^r\right)b}.
\end{equation*}
The above developments can be used to transform the eigenequation into the matrix eigenproblem 
\begin{equation} \label{eqmFICA}
\mathcal{G}^{1/2} \Sigma_{\tilde{A}\mathcal{G}^{1/2}}^{[4]} \mathcal{G}^{1/2} b=\hat{\kappa}_\theta\mathcal{G}_\theta b,
\end{equation}
where $\mathcal{G}_\theta=\mathcal{G}+\theta P^r$. Then, by performing the factorization $\mathcal{G}_\theta=LL^\top,$ the eigenequation  in Proposition  \ref{Prop02} can be rewritten as 
\begin{equation} \label{eqChol}
L^{-1}\mathcal{G}^{1/2} \Sigma_{\tilde{A}\mathcal{G}^{1/2}}^{[4]} \mathcal{G}^{1/2} (L^{-1})^\top v=\hat{\kappa}_\theta v,
\end{equation}
where $\Sigma_{\tilde{A}\mathcal{G}^{1/2}}^{[4]}$ is composed by \bl{$D_\theta=\ $diag$(\tilde{A}\tilde{\mathcal{G}}_\theta\tilde{A}^\top)$ with $\tilde{\mathcal{G}}_\theta=(L^{-1} \mathcal{G})^\top(L^{-1} \mathcal{G})$}, and $v=L^\top b,$ with $v^\top v =1.$

Now, defining $w= (L^{-1} \mathcal{G}^{1/2}) ^{-1} v,$ the eigenproblem turns on 
\begin{equation*}
\Sigma_{\tilde{A}\mathcal{G}^{1/2}}^{[4]} (L^{-1} \mathcal{G}^{1/2})^\top (L^{-1} \mathcal{G}^{1/2}) w = \hat{\kappa}_\theta w,
\end{equation*}
with $w^\top (L^{-1} \mathcal{G}^{1/2})^\top (L^{-1} \mathcal{G}^{1/2}) w =1 .$
This means that the smoothed functional ICA is equivalent to the ICA of matrix $A \mathcal{G}^{1/2}$ with a new metric in $\mathbb{R}^q$ defined by
$\left\langle x,y \right\rangle_{\mathcal{M}} = x^\top \mathcal{M} y$, for all $x,y \in \mathbb{R}^q.$

Therefore, solving (\ref{eqChol}) yields to $b_{\theta,j}=(L^{-1})^\top v_{j}$ such that $\hat{\psi}_{\theta,j}(t)=\phi(t)^\top b_{\theta,j}$ is the solution to the eigenequation \eqref{eqmFICA}. By computing the successive optimization problems in \eqref{1PSV} we obtain a set of orthonormal eigenfunctions verifying
\begin{equation*}
\Vert\hat{\psi}_{\theta,j}\Vert_{\theta}^{2}=b_{\theta,j}^{\top}
\mathcal{G}_{\theta}b_{ \theta,j}=v_{j}^{\top}v_{j} = 1
\end{equation*}
\begin{equation*}
\langle\hat{\psi}_{\theta,j},\hat{\psi}_{\theta,k}\rangle_{\theta}=b_{\theta,j}^{\top}\mathcal{G}_{\theta}b_{\theta,k}=v_{j}^{\top}v_{k}=0.
\end{equation*}
\end{proof}

\section{Interactive simulations} \label{secA1}
Interactive simulations regarding Sect.~\ref{TS4} are available through the Shiny app: \\ \href{ https://mvidal.shinyapps.io/fica/}{https://mvidal.shinyapps.io/fica/}.

\bibliographystyle{apalike-dashed}
 {\footnotesize
\bibliography{ref}}

\end{document}